%% file: submission-arXiv.tex
\documentclass[11pt]{amsart}
\usepackage[a4paper, margin=2.5cm]{geometry}
\usepackage{hyperref}
\hypersetup{ 
    colorlinks=true,
    linkcolor=blue,
    filecolor=magenta,      
    urlcolor=cyan,
}  
\usepackage{amssymb,bbm}
\def\fnsybm{2}

\usepackage{bbold}
\usepackage{xcolor}
\usepackage{cancel}
\usepackage{comment}
\includecomment{comment:techinical}
\makeatletter
\def\ifdraft{\ifdim\overfullrule>\z@
  \expandafter\@firstoftwo\else\expandafter\@secondoftwo\fi}
\makeatother
\newcommand{\draftcolor}[1]{\ifdraft{\color{#1}}{\color{#1}}}

\usepackage{tikz}
\usetikzlibrary{calc}

\newcommand{\mean}{\mathbb{E}}
\newcommand{\cL}{\mathcal{L}}
\newcommand{\bR}{\mathbb{R}}

\newcommand{\bone}{\mathbb{1}}
\newcommand{\bind}{\mathbf{1}}

\newcommand{\githubRemark}{The python algorithm to generate Figures \ref{fig:serv2}, \ref{fig:RENE.A2}, \ref{fig:RENE.C2}, \ref{fig:RENE.B2} and \ref{fig:RENE.D2} is avaiblabe for downloading at the public repository \cite{BDgithub}; see \cite{BDgithubpages} for an online implementation.}

\def\defeq{=}
\def\eqref#1{(\ref{#1})}
\def\state#1{\delta_{#1}}
\def\Kronecker#1#2{\{#1,#2\}}
\newcommand{\mD}{\Delta}
\newcommand{\tD}{{\tilde \mD}}
\newcommand{\mQ}{Q}
\newcommand{\tQ}{{\tilde \mQ}}
\newcommand{\diag}{{\rm diag}}
\newcommand{\phic}{{\phi_{*}}}
\newcommand{\psic}{{\psi_c}}
\newcommand{\tphic}{{{\tilde \phi}_{*}}}
\newcommand{\tpsic}{{{\tilde \psi}_c}}

\definecolor{candypink}{rgb}{0.89, 0.44, 0.48}
\definecolor{cadmiumgreen}{rgb}{0.0, 0.42, 0.24}

\newtheorem{theorem}{Theorem}
\newtheorem{lemma}{Lemma}
\newtheorem{corollary}{Corollary}
\newtheorem{remark}{Remark}
\newenvironment{proofof}[1]{\paragraph{\bf Proof of #1:}}{\qed}

\title{An $M/M/c$ queue with queueing-time dependent service rates}

\author{Bernardo D'Auria}
\address{Statistics Department, Madrid University Carlos III, Avda. Universidad 30, 28911 Leganés (Madrid), Spain}
\email[B.~D'Auria]{bernardo.dauria@uc3m.es}
\thanks{The research of the first author was partially supported by the Spanish Ministry of Economy and Competitiveness [Grant MTM2017-85618-P]}
\author{Ivo J.B.F. Adan}
\address{Industrial Engineering Department, Technische Universiteit Eindhoven, Postbus 513, 5600 MB, Eindhoven, The Netherlands}
\email[I.~Adan]{i.adan@tue.nl}
\author{René Bekker}
\address{Department of Mathematics, VU Amsterdam, De Boelelaan 1105, 1081 HV Amsterdam, The Netherlands}
\email[R.~Bekker]{r.bekker@vu.nl}
\author{Vidyadhar Kulkarni}
\address{Department of Statistics and Operations Research, University of North Carolina, Chapel Hill, NC 27599, USA}
\email[V.~Kulkarni]{vkulkarn@email.unc.edu}

\begin{document}
    \begin{abstract}
        Recent studies indicate that in many situations service times are affected by the experienced queueing delay of the particular customer.
        This effect has been detected in different areas,
        such as health care, call centers and telecommunication networks.
        In this paper we present a methodology to analyze a model having this property.
        The specific model is an $M/M/c$ queue in which any customer may be tagged at her arrival time if her queueing time will be above a certain fixed threshold. All tagged customers are then served at a given rate that may differ from the rate used for the non-tagged customers. 
        We show how it is possible to model the virtual queueing time of this queueing system by a specific Markov chain. Then, solving the corresponding balance equations, we give a recursive solution to compute the stationary distribution, which involves a mixture of exponential terms. Using numerical experiments, we demonstrate that the differences in service rates can have a crucial impact on queueing time performance.
    \end{abstract}
    \maketitle

\section{Introduction}
For classical queueing models, service times are typically assumed to be independent of experienced delay. Such independence assumptions are often crucial for analytical tractability of the queueing system's performance. 
In practice, however, it has been recognized that the amount of waiting affects service durations, and the assumed independence does therefore not hold. Empirical evidence of this dependence relation primarily stems from the health care domain. The studies \cite{batt2012, chalfin2007, chan2017, chan2008, renaud2009, richardson2002, siegmeth2005, soltani2019} indicate that delays in admission have adverse effects on patient outcomes and consequently increase the patients length of stay; this is referred to as the slowdown effect in \cite{selen2016}. 

At a more conceptual level, the field of behavioral operations investigates how servers and customers behave in an operational setting. 
The recent study \cite{delasay2019} indicates that in many situations service times are affected by the load. The authors develop a framework for the impact of load on service times, where they distinguish server, network, and customer mechanisms. For the server mechanism, it is observed (and supported by literature) that there is a clear impact of workload on the service speed, and this impact may go in different directions. The authors of \cite{delasay2019} found far fewer customer mechanisms in the literature, although they expect them to exist. A psychological view of a customers queueing experience during its sojourn time is provided in \cite{carmon1995}. Specifically, the authors assume that the dissatisfaction level of a customer increases during waiting, whereas this may be compensated during service. As a consequence, for an acceptable level of dissatisfaction after service, the service time should be longer for a customer experiencing longer delays. Moreover, after excessive waiting customers expect valuable service \cite{maister1984}, which may also affect the corresponding service time. Similarly, the recent study \cite{ulku2020} in a retail environment found that customers waiting longer in fact consume more. Thus, from the customer perspective, it seems conceivable that excessive waits are associated with longer service times. 

The aim of this paper is to find the steady-state queueing time distribution in multi-server queues where the service time is affected by the experienced queueing time. Despite its apparent practical relevance, such queues have hardly been studied in a service setting with multiple servers. More specifically, we let the service rate of each server depend on whether the experienced queueing time of the customer in service is above or below a given threshold upon service initiation. We envisage that this typically corresponds to a customer mechanism, although the service rate adaptation might also be the consequence of the server adapting to congestion. Despite the inherent model complexity, the steady-state queueing time distribution turns out to be remarkably tractable in this case
and can be expressed in terms of a mixture of exponential terms.
From our numerical experiments, we see that taking the differences in service rates into account results in crucially different queueing time behavior. As such, ignoring the dependence between experienced waiting and service time might be wholly inadequate. 
An online implementation of the model is available to further facilitate managerial decision making \cite{BDgithubpages}.

Delay thresholds are typically used in empirical health care studies to distinguish delayed and non-delayed patients; if the admission delay is above the delay threshold, a patient is considered to be delayed.
For instance, \cite{chalfin2007, richardson2002} investigated the impact of delayed patients at the emergency department on the inpatient length of stay. Based on the patient data, the difference in length of stay is in the order of hours. For cardiac patients, delays are even much more critical; delays in the order of minutes lead to adverse patient outcomes \cite{chan2008}. Less critical cases, such as surgery of hip fractures, have delay thresholds in the order of days \cite{siegmeth2005}. For patients with community-acquired pneumonia a similar delay threshold is used \cite{renaud2009}. For both situations it is shown again that delayed admissions experience extended length of stay. Another example of the impact of physician workload at the emergency department (ED) are \cite{batt2012, soltani2019}; amongst others, the authors observe that high physician workload leads to overtesting and generates extra post-ED care. 

The health care situations described above have recently inspired the study of multi-server queues, in which the service time (i.e. the length of stay) is affected by delay and congestion at the clinical ward, such as the Intensive Care Unit \cite{chan2017, dong2015, selen2016}. The study of \cite{chan2017} is also supported with data verifying the correlation between delay and length of stay. In \cite{chan2017}, the multi-server queue with delay-dependent service is abbreviated with $M/M(f)/c$; the focus from the queueing perspective is on approximations and bounds for the workload process. 
The multi-server variant with abandonments in the quality and efficiency driven (QED) regime is considered in \cite{dong2015}. Next to the fact that this involves an asymptotic analysis, the service rate adaptations are also instantaneous instead of the more intricate delay effects on individual customers. Such server mechanisms are referred to as operator slowdown in \cite{selen2016}, as opposed to customer slowdown. The model in \cite{selen2016} also involves a multi-server queue, where the service rate depends on whether a customer has to wait or not. In terms of the current paper, this means that the waiting threshold is at zero. In addition, \cite{selen2016} focuses on the number of customers instead of queueing times.

There have been some recent studies on multi-server queues where service times depend on delay. The authors of \cite{wu2019} consider a general multi-server queue with abandonments and derive fluid limits as a proxy for expected queueing times. Moreover, \cite{wu2019unpub} considers a setting with customer abandonments, where the service time is either endogenously or exogenously determined by the system's dynamics. The focus there is mainly on statistical estimation for both dependency situations. Finally, in \cite{do2018} the service speed is affected by behavioral factors, such as server speedup due to increased workload and social loafing when multiple workers share the workload. However, the analysis is in terms of queue lengths instead of queueing times. 

From the literature discussed above, we observe that almost all studies of multi-server queues with delay-dependent service involve some sort of approximation. This is different for the single-server case, which is much more amenable for analysis. An important observation for the single-server case is that the queueing time then corresponds to the workload a customer finds upon arrival; this is no longer the case for the multi-server setting with delay-dependent service.
There is a long tradition of single-server queues with workload-dependent features; we refer to \cite{dshalalow1997} for an early overview containing many references. Among those early papers are \cite{posner1973} and \cite{brill1981}. Interestingly, in 1973 Posner already noted that the server may provide more appropriate service to counter the negative effect of waiting \cite{posner1973}; the author then provides a complete analysis for the M/M/1 case in which the service rate is a step function of the queueing time. A little later, \cite{brill1981} provides an exact analysis for the M/M/2 queue where non-waiting customers have a different service rate. 

For workload-dependent M/G/1 queues, often the service and/or arrival rates are assumed to depend on the workload, but not so often the complete service time. However, generalizations of such systems are L\'evy driven queues in which the L\'evy exponent depends on the position of the process. The L\'evy exponent incorporates the Laplace transform of the service time distribution and, hence, the service time may thus depend on the workload found by a customer entering service. Examples of such L\'evy driven queues with state-dependent exponent are \cite{bekker2009, bekker2009levy, palmowski2011}. 
Finally, \cite{whitt1990} and \cite{boxma2007} consider G/G/1 queues with service and interarrival times that depend linearly on delays.

Limiting distributions in terms of mixtures of exponentials are also common in Markov-modulated fluid models. In fact, our analysis is along similar lines as such fluid models, although our differential equations differ from the ones found in traditional fluid queues \cite{anick1982}, see \cite{kulkarni1997} for an early overview. Some examples of fluid models with level-dependent features are \cite{dasilva2009, malhotra2009, scheinhardt2005}. 
A crucial difference with fluid models is the role of the background state. 
Our state description, where the service time depends on experienced delay, is delicate. In our case, the background state should be interpreted as the server state process; our state description is based on \cite{adan2019}. 

The paper is organized as follows. In Section~\ref{sec:model}, a model and state description is provided. Section~\ref{sec:balance} presents balance equations that are required to determine the limiting distribution. The limiting distribution is derived in Section~\ref{sec:lim}, followed by some illustrative special cases in Section~\ref{sec:cases}. Section~\ref{sec:numerics} contains some numerical insights. For readability, most of the technical proofs are deferred to Appendix~\ref{app:proofs}. 
A python algorithm to compute the queueing time distribution is avaiblabe for downloading at the public repository \cite{BDgithub}; see \cite{BDgithubpages} for an online implementation.

\section{Model and state description} \label{sec:model}
We consider a queueing system with $c$ identical servers and an infinite waiting room. Customers arrive according to a Poisson process with rate $\lambda$. Let $W(t)$ be the virtual queueing time (VQT) at time $t$. That is, if a customer arrives at time $t$, his service will start at time $t+W(t)$. Clearly, if at least one server is idle at time $t$, $W(t)=0$. If all servers are busy at time $t$, $W(t) > 0$. The service times of the customers depend on their queueing time through a critical level $k > 0$ as follows: if a customer arrives at time $t$, and $W(t) \le k$, he is classified as a class 1 customer, and his service time is exp($\mu_1$), otherwise he is classified as class 2 customer and his service time is exp($\mu_2$). In order to describe the dynamics of the VQT process $\{W(t), t \ge 0\}$, we introduce the server state process $S(t) = (S_1(t),S_2(t))$ as follows. We say that $S(t) = (i,j)$ if $i$ servers are serving class 1 customers and $j$ servers are serving class 2 customers at time $t+W(t)$, just before the new service starts at time $t+W(t)$. Clearly, we must have $0 \le S_1(t) + S_2(t) \le c-1$ for all $t \ge 0$. Furthermore,
\[W(t) > 0 \Rightarrow S_1(t) + S_2(t) = c-1.\]
and
\[0 \le S_1(t) + S_2(t) < c-1 \Rightarrow W(t) = 0.\]
We discuss the evolution of the $\{(W(t),S_1(t),S_2(t)), t \ge 0\}$ process below. We will use the following notation for the aggregate service rate:
\begin{equation}\label{def:Delta.func}
\Delta(i,j) = i\mu_1+j\mu_2 .
\end{equation}

Suppose the state at time 0 is $(0,i,j)$ with $0 \le i+j < c-1$. If the next event is an arrival, the state jumps to $(0,i+1,j)$; if it is  a departure of type 1, it jumps to state $(0,i-1,j)$; and if it is a departure of type 2, it jumps to state $(0,i,j-1)$. Hence, the transition rate from state $(0,i,j)$ to state $(0,i+1,j)$ is $\lambda$, to state $(0,i-1,j)$ is $i\mu_1$ and to state $(0,i,j-1)$ is $j\mu_2$.

Next, suppose the state at time 0 is $(0,i,c-1-i)$ with $0 \le i \leq c-1$. Again, if the next event is a departure of type 1, it jumps to state $(0,i-1,c-1-i)$; and if it is a departure of type 2, it jumps to state $(0,i,c-1-2)$. Hence, the transition rate from state $(0,i,c-1-i)$ to state $(0,i-1,c-1-i)$ is $i\mu_1$ and to state $(0,i,c-i-2)$ is $(c-1-i)\mu_2$. If the next event is an arrival, all servers become busy, $i+1$ of them serving type 1 customers and $(c-1-i)$ of them serving type 2 customers. The next departure occurs after an exp($\Delta(i+1,c-1-i)$) amount of time and the VQT process jumps to level $X \sim $ exp($\Delta(i+1,c-1-i)$). Also, the next departure is of type 1 with probability $(i+1)\mu_1/\Delta(i+1,c-1-i)$ and of type 2 with probability $(c-1-i)\mu_2/\Delta(i+1,c-1-i)$. Combining these observations, we see that the transition rate to state $(x,i,c-1-i)$ is
\[ \lambda (i+1)\mu_1 \exp(-\Delta(i+1,c-1-i)x)dx,\]
and to state  $(x,i+1,c-i-2)$ is 
\[ \lambda (c-1-i)\mu_2 \exp(-\Delta(i+1,c-1-i)x)dx.\]
This completes the description of all transitions out of states $(0,i,j)$ with $0 \le i+j \le c-1$. 

Next, consider states $(w,i,c-1-i)$ with $0 < w \le k$ and $0 \le i \le c-1$. The state does not change if the next event is a departure. It can change only if the next event is an arrival. An arrival in this state is of type 1.  By following the same argument as in the case of state $(0,i,c-1-i)$, we see that the transition rate to state $(w+x,i,c-1-i)$ is
\[ (i+1)\mu_1 \exp(-\Delta(i+1,c-1-i)x)dx,\]
and  to state $(w+x,i+1,c-i-2)$ is 
\[ (c-1-i)\mu_2 \exp(-\Delta(i+1,c-1-i)x)dx.\]

Now consider states $(w,i,c-1-i)$ with $w > k$ and $0 \le i \le c-1$. An arrival in this state is of type 2. Hence, following the same argument as above, we see that the transition rate to state $(w+x,i,c-1-i)$ is 
\[ (c-i)\mu_2 \exp(-\Delta(i,c-i)x)dx,\]
and  to state $(w+x,i-1,c-i)$ is 
\[ i\mu_1 \exp(-\Delta(i,c-i)x)dx.\]

Finally, if $W(t) > 0$, the VQT process changes continuously at rate $-1$ between arrivals. This completes the description of the evolution of the  process $\{(W(t), S(t)), t \ge 0\}$.

\section{Limiting distribution: balance equations}
\label{sec:balance}

In this section, we derive the balance equations for the VQT process $(W(t), S(t))$ defined in Section~\ref{sec:model} that are satisfied by the limiting distribution.
It is straightforward to see that the VQT process is stable if 
\begin{equation}\label{eq:stability.condition} \frac{\lambda}{\mu_2 c} < 1.
\end{equation}
We shall assume stability from now on and focus on the limiting distribution of $(W(t), S(t))$.\\

Now let, for $x \ge 0$ and $t \ge 0$,
\[ F_i(t,x) = P(0 < W(t) \le x; S(t)=(i,c-1-i)), \;\;\; 0 \le i \le c-1. \]
Define $F_i(x) = \lim_{t \rightarrow \infty} F_i(t,x)$, and define the row vector function
\[ F(x) = [F_0(x), \; F_1(x), \; \cdots, F_{c-1}(x)],\]
whose first two derivatives are denoted by $F'(x)$ and $F''(x)$.
Also, for the case that no customers are waiting, let
\[ \pi(i,j) = \lim_{t \rightarrow \infty} P(W(t)=0, S(t) = (i,j)),\;\;\; 0 \le i+j \le c-1 , \]
and
\[ \state{i} = [\pi(j,i-j), \, 0 \le j \le i], \;0 \le i \le c-1 . \]

Using the transition rates derived above, we see that the  $\pi$'s satisfy the following balance equations:
\begin{eqnarray}
(\lambda + i\mu_1 + j\mu_2) \pi(i,j) & = & (i+1)\mu_1 \pi(i+1,j) + (j+1)\mu_2 \pi(i,j+1)  \nonumber\\ 
& & + \lambda \bind_{(i>0)} \pi(i-1,j), \;\;\; 0 \le i+j < c-1. \label{eq:piij0}
\end{eqnarray}
with $\bind_{(\cdot)}$ denoting the indicator function. 

Next we derive the integro-differential equations satisfied by $F(\cdot)$ for the case there is queueing delay. 

We denote by $I$ the identity matrix, whose size will be clear from the context, 
and by $\hat I$ a rectangular matrix obtained from $I$ by adding a null column on the left, i.e. $\hat I = (\begin{array}{ll} 0 & I \end{array})$.
Throughout the paper, we will use the convention of denoting by $\hat \cdot$ a rectangular  (instead of square) matrix.

Let $B_1$ be a $c \times c$ square matrix with entries given by
\begin{eqnarray*}
B_1(i,i) & = & (i+1)\mu_1, \;\; 0 \le i \le c-1, \\
B_1(i,i+1) & = & (c-i-1)\mu_2, \;\;\; 0 \le i < c-1,\\
 B_1(i,j)& = & 0 \;\;\;\mbox{for all other } (i,j),
\end{eqnarray*}
and $B_2$ be a $c \times c$ square matrix with entries given by
\begin{eqnarray*}
B_2(i,i) & = & (c-i)\mu_{2}, \;\; 0 \le i \le c-1, \\
B_2(i,i-1) & = & i\mu_{1}, \;\;\; 1 \le i \le c-1,\\
B_2(i,j) & = & 0\;\;\; \mbox{for all other } (i,j).
\end{eqnarray*}

We finally define the matrices 
\begin{equation}\label{def:mD.i}
\mD_i =  \mbox{diag}(j\mu_1 + (i-j)\mu_2, \, 0 \le j \le i) , \;\; 0 \le i \le c-1, 
\end{equation}
and 
\begin{eqnarray}
\tD_\kappa & = & \mu_\kappa I + B_\kappa^{-1}(\mD_{c-1}) B_\kappa , \;\; \kappa\in\{1,2\}, \label{def:tD.kappa} \\
\tQ_\kappa(x) & = & \exp(-\tD_\kappa x) , \;\; \kappa\in\{1,2\}.  \label{def:tQ.kappa} 
\end{eqnarray}

\begin{theorem} \label{th:dFx}
The limiting distribution vector $F$ satisfies the following integro-differential equations:
\begin{eqnarray}
F'(x) & = & 
\lambda F(x)  - \lambda \int_0^x F(y) B_1 \tQ_1(x-y) dy + F'(0) \nonumber\\
& & - \lambda \state{c-1} B_1 (I - \tQ_1(x)) \tD_1^{-1} ,\;\;\; 0 < x < k  \ , \label{eq:Q1}\\
F'(x) & = &  \lambda F(x)    
   - \lambda \int_k^x F(y) B_2 \tQ_2(x-y) dy + F'(0) \nonumber\\
& & -  \lambda \state{c-1} B_1 (I - \tQ_1(x)) \tD_1^{-1}
     - \lambda \int_0^k F(y) B_1 \tQ_1(x-y)  dy \nonumber\\
& & - \lambda F(k) \int_k^x (B_1 \tQ_1(x-y)  - B_2 \tQ_2(x-y) ) dy
,\;\;\; x > k \ . \label{eq:Q2}
\end{eqnarray}
\end{theorem} 

\begin{proof}
    The proof follows standard probabilistic arguments 
    and uses an infinitesimal approach, which we defer
    to the \hyperref[pf:th:dFx]{Appendix}.
\end{proof}

Equations~\eqref{eq:Q1} and \eqref{eq:Q2} are integro-differential equations. 
These equations are related to level crossings principles; see Subsection~\ref{subs:single} for the single-server case providing additional intuitive insight.
To find the limiting distribution, we first convert them to second order linear non-homogeneous differential equations. They are given in the following theorem. To do so, first define the differential operators $\cL_1$ and $\cL_2$ as follows:
\[\cL_\kappa G(x) 
    = G''(x)  
    - G'(x) (\lambda I - \tD_\kappa) 
    + \lambda G(x) (B_\kappa - \tD_\kappa) \quad \kappa\in\{1,2\}.\]

\begin{theorem} \label{th:dFx2}
The limiting distribution vector $F$ satisfies the following second order differential equations:
    \begin{eqnarray}
    \cL_1 F(x)
    &=& \alpha_0
    ,\;\;\; 0 < x < k  \ , \label{eq:F12} \\
    \cL_2 F(x)
    &=& \alpha_1 + \alpha_2 \tQ_1(x-k)(\tD_1 - \tD_2) 
    ,\;\;\; x > k \ , \label{eq:F22} 
    \end{eqnarray}
    where 
    \begin{eqnarray}
        \alpha_0 & = &  F'(0) \tD_1 - \state{c-1} \lambda B_1  \ , \label{eq:alp0} \\
        \alpha_1 & = &  \alpha_0 \tD_1^{-1} \tD_2 - \lambda F(k+) (B_1\tD_1^{-1}\tD_2  - B_2) \ , \label{eq:alp1}\\
        \alpha_2 & = &  \alpha_1 \tD_2^{-1} - F'(k+) + \lambda F(k+) (I - B_2\tD_2^{-1}).\label{eq:alp2}
    \end{eqnarray}
\end{theorem}

\begin{proof} 
This follows from rewriting the integro-differential equations \eqref{eq:Q1} and  \eqref{eq:Q2}, see \hyperref[pf:th:dFx2]{Appendix}.
\end{proof}

The next corollary gives the boundary conditions for $F$. Here and later we use the notation $g(x\pm):=\lim_{y \to x^\pm} g(y)$ to denote the one-sided limits of the function $g$ at $x$.

\begin{corollary}\label{cor:bc}
    The limiting distribution vector $F$ satisfies the following boundary conditions:
    \begin{eqnarray}
    F(0) & = & 0, \label{Fzero}\\
        F(k-) & = & F(k+), \label{F.cont.cond}\\
        F'(k-) & = & F'(k+),  \label{F.diff.cond}\\
         F'(0) & = & \state{c-1}(\lambda I + \mD_{c-1}) - \state{c-2} \lambda \hat I. \label{F.zero.cond}
    \end{eqnarray}
\end{corollary} 

\begin{proof}
    See \hyperref[pf:cor:bc]{Appendix} for details.
\end{proof}

\section{Solution of the balance equations} \label{sec:lim}
In this section we determine the limiting distribution by 
developing the analytical solution of the differential equations in Theorem~\ref{th:dFx2} 
and the boundary conditions in Corollary~\ref{cor:bc}.
In fact, we present two different ways to express the limiting distribution of the VQT process. The first is based on a scalar representation and clearly reveals that $F(x)$ can be written as a mixture of exponentials. This representation gives insight in the probabilistic interpretation of the queueing delay and is presented in Subsection~\ref{subs:lim}. 
The second concerns a matrix representation and is more compact. This representation is more amenable for numerical computations and can be found in Subsection~\ref{subs:comp}. \\

\subsection{Limiting distribution} \label{subs:lim}
For the solution of Equation~(\ref{eq:F12}) in Theorem~\ref{th:dFx2}, we first need to solve the homogeneous equation. Hence, we first need to define several backgrounds quantities. Consider the  following  quadratic eigenvalue equation: 
\begin{equation}\label{eq:eig1}
\phi[\theta^2 I - \theta (\lambda I - \tD_1) + \lambda (B_1 - \tD_1) ] = 0.
\end{equation}
There are $2c$ solutions $\{(\theta_i,\phi_i), \; 0 \le i \le 2c-1\}$ to the above system. Since the matrices involved in the above equation are all upper triangular, it is easy to see that these $2c$ solutions are given by the solutions to the following quadratic equations:
\begin{equation}\label{eq:theta}
\theta^2 - \theta(\lambda - (i+1)\mu_1 - (c-1-i)\mu_2) -(c-1-i)\lambda\mu_2 = 0, \;\;\; 0 \le i \le c-1.
\end{equation}
If both conditions $\lambda \not= c\mu_1$ and $\lambda \not= c(\mu_1 - \mu_2)$ are satisfied, all these eigenvalues are real and distinct, see also Remark~\ref{rm:identicalEV} for the other cases.
For the rest of the paper we implicitly assume that these conditions hold. The solutions $\{\theta_i, 0 \le i \le c-1\}$ are given by
\begin{eqnarray}\label{eq.theta.neg}
    \lefteqn{\theta_i = \frac{1}{2}\left(\lambda - (i+1)\mu_1 - (c-1-i)\mu_2\right)} \nonumber \\
    & & - \frac{1}{2}\sqrt{(\lambda - (i+1)\mu_1 - (c-1-i)\mu_2)^2 +4(c-1-i)\lambda\mu_2}. 
\end{eqnarray}
Note that $\{\theta_i, 0 \le i \le c-2\}$ are negative and 
$\theta_{c-1} = \min\{0, \lambda - c\mu_1\}$.
The solutions $\{\theta_{i+c}, 0 \le i \le c-1\}$ are positive and are given by
\begin{eqnarray}\label{eq.theta.pos}
    \lefteqn{\theta_{i+c} = \frac{1}{2}\left(\lambda - (i+1)\mu_1 - (c-1-i)\mu_2\right)} \nonumber \\
    & & + \frac{1}{2}\sqrt{(\lambda - (i+1)\mu_1 - (c-1-i)\mu_2)^2 +4(c-1-i)\lambda\mu_2}.
\end{eqnarray}
Note that $\theta_{2c-1} = \max\{0, \lambda - c\mu_1\}$. 
The corresponding eigenvectors $\{\phi_i, 0 \le i \le 2c-1\}$ are easy to compute. 
In particular, the eigenvector corresponding to the null eigenvalue is denoted by  
 \begin{equation} \label{eq:phic} 
    \phic = [0, \; 0, \cdots, 0, \; 1],
 \end{equation}
 which is a row vector of length $c$.

Next, we turn to the homogeneous equation based on Equation~(\ref{eq:F22}). For this, consider the following quadratic eigenvalue equation:
\begin{equation} \label{eq:beta} 
    \psi[\beta^2 I - \beta (\lambda I - \tD_2) + \lambda (B_2 - \tD_2) ] = 0.
\end{equation}
There are $2c$ solutions $\{(\beta_i,\psi_i), \; 0 \le i \le 2c-1\}$ to the above system. The $\beta_i$'s for $0 \le i \le c-1$ are given by 
\begin{equation} \label{eq.beta.neg}
    \beta_i = \frac{1}{2}\left(\lambda - i\mu_1 - (c-i)\mu_2 - \sqrt{(\lambda - i\mu_1 - (c-i)\mu_2)^2 +4i\lambda\mu_1}\right),
\end{equation}
and the $\beta_{i+c}$'s, for $0 \le i \le c-1$ are given by 
\begin{equation} \label{eq.beta.pos}
    \beta_{i+c} = \frac{1}{2}\left(\lambda - i\mu_1 - (c-i)\mu_2 + \sqrt{(\lambda - i\mu_1 - (c-i)\mu_2)^2 +4i\lambda\mu_1}\right).
\end{equation}
Similarly to \eqref{eq.theta.neg}, and \eqref{eq.theta.pos}, all these eigenvalues are real and distinct, assuming that $\lambda \not= c(\mu_2 - \mu_1)$. The eigenvalues $\{\beta_i, 0 \le i \le c-1\}$ are negative,
$\beta_c = 0$ and $\{\beta_{i+c}, 1 \le i \le c-1\}$ are positive. 
The corresponding eigenvectors $\{\psi_i, 0 \le i \le 2c-1\}$ are easy to compute. In particular, 
 \begin{equation} \label{eq:psic} 
    \psic = [1, \; 0, \cdots, 0, \; 0],
 \end{equation}
 which is a row vector of length $c$.\\

\begin{remark} \label{rm:eigenvalues}
The parameters $\theta_i$ and $\beta_i$ may be related to the virtual waiting in regular M/M/$c$ queues. For instance, $\rm{exp}(\theta_{c-1}) = \rm{exp}(\lambda - c\mu_1)$ and $\rm{exp}(\beta_{0}) = \rm{exp}(\lambda - c\mu_2)$ are proportional to the stationary densities of the VQT in M/M/$c$ queues with only service rates $\mu_1$ and $\mu_2$, respectively.
Moreover, consider the process $W(t) \in (0,k)$ and fix the server state process $S(t)=(i,c-1-i)$; then the VQT process is decreasing with rate 1 and makes jumps with rate $\lambda$ of size exp($\Delta(i+1,c-1-i)$). Upon a jump, the server state process $S(t)$ changes with probability $(c-1-i)\mu_2/\Delta(i+1,c-1-i)$, which may be interpreted as a type of clearing \cite{boxma2001}. It may be verified that the stationary density of such a `clearing system' is a mixture of exp($\theta_i$) and exp($\theta_{i+c}$). 
A similar argument applies to $W(t) > k$ in terms of $\beta_i$.
\end{remark}

As mentioned, to get the solution of the differential equations of Theorem~\ref{th:dFx2} we first need to find the  solution of the homogeneous differential equations using the above, and then we look for a particular solution.
%
    However, in order to construct a particular solution, since both equations 
    \eqref{eq:F12} and \eqref{eq:F22} admit zero as eigenvalue, we would need together with the left eigenvectors
    $\phic$ and  $\psic$, respectively defined in \eqref{eq:phic} and \eqref{eq:psic}, 
    the corresponding right eigenvectors that we denote by $\tphic$ and $\tpsic$.
The following result shows an important relation between those eigenvectors that will be used later in the proof of Theorem~\ref{th:Fx2sol} to show that $F(x)$ does not have a linear term.
    
    \begin{lemma}\label{lm:right.eigenvec}
        Let $\tphic$ and $\tpsic$ be the right eigenvectors corresponding to the left eigenvectors $\phic$ and $\psic$, i.e. satisfing the following relations
        \begin{eqnarray} 
            (B_1 - \tD_1) \tphic & = & 0, \label{eq:tphic.prop} \\
            (B_2 - \tD_2) \tpsic & = & 0. \label{eq:tpsic.prop} 
        \end{eqnarray}
        It follows that 
        \begin{equation}\label{eq:right.eigenvec}
            \alpha_1 \cdot \tpsic = 0 \Rightarrow \alpha_0 \cdot \tphic = 0 .
        \end{equation}
    \end{lemma}   
    \begin{proof} 
        The proof uses linear algebra techniques and is included in the \hyperref[pf:lm:right.eigenvec]{Appendix}.
    \end{proof} 

    The next result gives the solution of the differential equations of Theorem~\ref{th:dFx2} 
    in terms of 4$c$ unknowns $\{a_i, 0 \le i \le 2c-1\}$ and $\{b_i, 0 \le i \le c \}$. 

\begin{theorem} \label{th:Fx2sol}
For $0 < x < k$, the solution is given by
\begin{equation} \label{eq:gsol1}
F(x) = \sum_{i=0}^{2c-1}a_ie^{\theta_ix}\phi_i + \alpha_0 M_0, \;\; 0 \le x \le k,
\end{equation}
with
\begin{equation} \label{eq:M0}
    M_0 = (\lambda (B_1 - \tD_1) + \diag(\phic))^{-1}. 
\end{equation}
For $x > k$, the solution is given by
\begin{equation} \label{eq:gsol2}
F(x) = \sum_{i=0}^{c-1}b_ie^{\beta_i (x-k)}\psi_i + b_c \psic 
+ \alpha_1 M_1 + \alpha_2 \tQ_1(x-k) (\tD_1-\tD_2) M_2, 
\end{equation}
with
\begin{eqnarray} 
    M_1 &=& (\lambda (B_2 - \tD_2) + \diag(\psic))^{-1},  \label{eq:M1} \\
    M_2 &=& ((c \mu_1 + \lambda)(c \mu_1 I-\tD_2)+\lambda B_2)^{-1}. \label{eq:M2} 
\end{eqnarray}
\end{theorem}

\begin{proof} 
This follows from solving the systems of second order linear differential equations below and above level $k$ in Theorem~\ref{th:dFx2}, thereby also utilizing Lemma~\ref{lm:right.eigenvec}; see \hyperref[pf:th:Fx2sol]{Appendix}.
\end{proof}

\begin{remark} \label{rm:identicalEV}
We note that in the special case that $\lambda = c\mu_1$, we have two identical eigenvalues $\theta_{c-1}=\theta_{2c-1}=0$. 
Furthermore, in case $\lambda = c(\mu_1 - \mu_2)$, it holds that $\theta_i = \theta_{c-1} = \lambda - c\mu_1$ for all $i=0,\ldots,c-1$. In that case, $F(x)$ for $0<x<k$, contains terms of the form $\sum_{i=0}^{c-1} a_i x^i e^{(\lambda - c\mu_1)x}$. Similarly, if $\lambda = c(\mu_2 - \mu_1)$, then the $\beta_i$ for $i=0,\ldots,c-1$ are identical to $\beta_0$ and the mixture of exponentials in $F(x)$ for $x>k$ needs to be replaced. 
\end{remark}

\begin{remark} \label{rm:mu1=mu2}
In case $\mu_1=\mu_2$ the model reduces to a simple M/M/c queue. Therefore, in this case the VQT distribution is given by 
        $$P(W \le x) = 1 - C(c, \lambda/\mu_2) \frac{1}{1-\rho} \exp\{- c\mu_2 (1 - \rho) x\} , $$
where $C(c, \lambda/\mu_2) = 1/\left(1+(1-\rho)(\frac{c!}{c\rho^c})\sum_{k=0}^{c-1}\frac{(c\rho)^k}{k!}\right)$ is the Erlang's C formula
and $\rho=\lambda/(c\mu_2)$ is the \emph{server utilization}. 
\end{remark}

We next use the boundary conditions in Corollary~\ref{cor:bc} to solve for the $3c+1$ unknown constants $\{a_i, 0 \le i \le 2c-1\}$ and $\{b_i, 0 \le i \le c\}$. We also have $c(c+1)/2$ probabilities $\pi(i,j), \; 0 \le i+j \le  c-1$ that need to be determined. The result is given in the next theorem.

\begin{theorem} \label{th:ab}
The constants $\{a_i, 0 \le i \le 2c-1\}$, $\{b_i, 0 \le i \le c\}$ and  probabilities $\pi(i,j), \; 0 \le i+j \le c-1$ satisfy the following equations:
\begin{equation} \label{eq:con1} \sum_{i=0}^{2c-1}a_i\phi_i + \alpha_0 M_0 = 0,
\end{equation}
\begin{equation} \label{eq:con2} \sum_{i=0}^{2c-1}a_ie^{\theta_i k}\phi_i + \alpha_0 M_0 = \sum_{i=0}^{c-1}b_i\psi_i + b_c \psic 
+ \alpha_1 M_1 + \alpha_2 (\tD_1-\tD_2) M_2,
\end{equation}
\begin{equation} \label{eq:con3}\sum_{i=0}^{2c-1}a_i\theta_ie^{\theta_i k}\phi_i = \sum_{i=0}^{c-1}b_i\beta_i\psi_i - \alpha_2 \tD_1 (\tD_1-\tD_2) M_2,
\end{equation}
\begin{equation} \label{eq:con4} \sum_{i=0}^{2c-1}a_i\theta_i\phi_i = \state{c-1}(\lambda I + \mD_{c-1}) - \state{c-2} \lambda \hat I, 
\end{equation}
\begin{eqnarray} 
(\lambda + i\mu_1 + j\mu_2) \pi(i,j) & = & (i+1)\mu_1 \pi(i+1,j) + (j+1)\mu_2 \pi(i,j+1)  \nonumber\\ 
& & + \lambda \bind_{(i>0)} \pi(i-1,j), \;\;\; 0 \le i+j < c-1 \label{eq:con5}
\end{eqnarray}
and the normalizing equation
\begin{equation} \label{eq:con6}b_c \psic \bone + \alpha_1 M_1 \bone + \sum_{i,j:i+j=0}^{i+j=c-1}\pi(i,j) = 1,
\end{equation}
where $\bone$ denotes the all-one vector.
\end{theorem}

\begin{proof} Equations~\eqref{eq:con1}--\eqref{eq:con4} follow from the boundary conditions~\eqref{Fzero}--\eqref{F.zero.cond}, respectively.
Equation~\eqref{eq:con5} presents the balance equations in \eqref{eq:piij0}.
Finally, Equation~\eqref{eq:con6} is the normalizing equation. 
\end{proof}

Equations~\eqref{eq:con1} to \eqref{eq:con6} from the above theorem yield $4c + c(c-1)/2 + 1 = 3c+1+c(c+1)/2$ linear equations for the  $3c+1$ unknown constants $\{a_i, 0 \le i \le 2c-1\}$, $\{b_i, 0 \le i \le c\}$ and the $c(c+1)/2$ unknown  probabilities $\pi(i,j), \; 0 \le i+j \le c-1$.
Theorem \ref{th:ab} gives necessary conditions for the constants to satisfy, but it does not assure that Equations (\ref{eq:con1}--\ref{eq:con6}) characterize them. In the following section we show, starting from these equations, how we are able to construct the unique solution.

\subsection{Computing the solution} \label{subs:comp}
For the limiting distribution in Theorem~\ref{th:Fx2sol}, we need 
Theorem \ref{th:ab} that expresses all constants as a solution of a quite large system of linear equations.
In this section we try to express the solution in simpler terms that are easier to implement in a computer language equipped with basic matrix functions.
In addition the contruction of the solution given  below shows that the system in Theorem \ref{th:ab} only admits a unique solution.

First of all, it helps to write the function $F(x)$, given in \eqref{eq:gsol1} and \eqref{eq:gsol2},
in matrix form.
We start by defining the matrices $\Phi^- = [\phi_i^t]^t$, $0 \le i \le c-1$ and $\Phi^+ = [\phi_i^t]^t$, $c \le i \le 2c-1$, whose rows are given by the eigenvectors corresponding to the eigenvalues in \eqref{eq.theta.neg} and \eqref{eq.theta.pos}. We also define the diagonal matrices $\Theta^- = \mbox{diag}(\theta_i)$, $0 \le i \le c-1$, and $\Theta^+ = \mbox{diag}(\theta_i)$, $c \le i \le 2c-1$.
Then, the matrices 
        $U^\pm_1 = (\Phi^\pm)^{-1} \Theta^\pm \Phi^\pm$,
        solve the equation
        \[ U_1^2 - U_1 (\lambda I - \tD_1) + \lambda (B_1 - \tD_1)  = 0, \]
        with $U^-_1$~having all non-positive (negative) eigenvalues, $U^+_1$~having all positive (non-negative) eigenvalues if $\lambda>c\mu_1$ ($\lambda<c\mu_1$).

        In a similar way we construct the matrices $U^\pm_2$ solving the equation 
        \[ U_2^2 - U_2 (\lambda I - \tD_2) + \lambda (B_2 - \tD_2)  = 0 \]
        with $U^-_2$~with all negative eigenvalues and $U^+_2$~with all non-negative eigenvalues.
        This allows us to rewrite the expression in \eqref{eq:gsol1} and \eqref{eq:gsol2} as
        \begin{eqnarray*}
        F(x) &=& a^- e^{U^-_1 x} + a^+ e^{U^+_1 x} + \alpha_0 M_0 
            , \;\;\; 0 \le x \le k, \\
        F(x) &=& b^- e^{U^-_2 x} + b_c \psic 
            + \alpha_1 M_1 + \alpha_2 \tQ_1(x-k) (\tD_1-\tD_2) M_2
            , \;\;\; x \ge k ,
        \end{eqnarray*}
        for unknown constant vectors $a^-$, $a^+$ and $b^-$.
        We then express these vectors in terms of the unknown vector $F'(0)$ by using the continuity conditions given in Corollary~\ref{cor:bc}.
        This gives the following easier expressions
\begin{eqnarray}\label{eq:F1.exp}
    F(x) &=& (F'(0) + \alpha_0 M_0 U^-_1) (U^+_1 - U^-_1)^{-1} (e^{U^+_1 x} - e^{U^-_1 x}) \nonumber \\
    & & + \alpha_0 M_0 (I - e^{U^-_1 x}), \;\;\; 0 \le x \le k , \\
    \label{eq:F2.exp}
    F(x) &=& F(k-) e^{U^-_2 (x-k)} + (b_c \psic + \alpha_1 M_1) (I-e^{U^-_2 (x-k)})  \nonumber \\
    & & - \alpha_2 (\tD_1-\tD_2) M_2 e^{U^-_2 (x-k)} \nonumber \\
    & & + \alpha_2 \tQ_1(x-k) (\tD_1-\tD_2) M_2, \;\;  x \ge k .
\end{eqnarray}
In the theorem below, $F'(0)$ is expressed in terms of $\state{c-1}$ and $b_c$.

\begin{theorem} \label{th:FxMsol}
    The function $F$, given in \eqref{eq:F1.exp} and \eqref{eq:F2.exp}, may be written in terms of the constant $b_c$ and the vector $\state{c-1}$, since
    \begin{equation}
        F'(0) = \state{c-1} H_{16} - b_c \psic   H_{15}, \label{eq:F'0}
    \end{equation}
    where $H_{15}$ and $H_{16}$ are defined in \eqref{eq:H15} and \eqref{eq:H16}, respectively.
    
    In particular it follows that 
    \begin{equation}
        F(\infty) = \state{c-1} H_{20} + b_c \psic  H_{19},  \label{eq:Finfty}
    \end{equation}
    with $H_{19}$ and $H_{20}$ defined in \eqref{eq:H19} and \eqref{eq:H20}, respectively.
\end{theorem}

\begin{proof} 
    This follows from Corollary~\ref{cor:bc} and some tedious rewriting, 
    see \hyperref[pf:th:FxMsol]{Appendix}.
\end{proof}

Having expressed the function $F(x)$ for $x\ge0$, in terms of $\state{c-1}$, 
it is only left to find the probability of the discrete states, together with the constant $b_c$ that can be found by using the normalizing equation.
This is the result of following theorem.

\begin{theorem} \label{th:bc.sol}
    The discrete probabilities can be computed as follows
    \begin{equation}
        \label{eq:state.sol} 
        \state{i}  = b_c \psic \hat H_{i} ,  \;\;\;  0 \leq i \leq c-1,
    \end{equation}
    where the matrices $\hat H_{i}$ are defined in \eqref{eq:hat.H} 
    and with the constant $b_c$ computed as
\begin{equation}
    \label{eq:bc.sol} 
    b_c^{-1} = \psic \left( H_{19} + \hat H_{c-1} H_{20} + \sum_{0 \le n \le c-1}  \hat H_n\right) \, \bone ,
\end{equation}
where $H_{19}$ and $H_{20}$ are given in \eqref{eq:H19} and \eqref{eq:H20} in the Appendix.
\end{theorem}

\begin{proof} 
The linear system of equations may be rewritten to recursively express $\state{c-1}$ in terms of $b_c$, whereas $b_c$ follows from normalization;   
see \hyperref[pf:th:bc.sol]{Appendix}.
\end{proof}

\section{Special cases}
\label{sec:cases}
In this section, we present two special cases that provide probabilistic understanding of the limiting VQT process. In Subsection~\ref{subs:single} we focus on the single-server case, whereas Subsection~\ref{subs:c2} covers the case with 2 servers.

\subsection{Single-server queue} \label{subs:single}
\begin{figure}
\centering
\begin{tikzpicture}[yscale=0.6, xscale=0.7]
    \coordinate (O) at (0,0);
    \draw [<->,thick] (0,8) node (w-axis) [above] {$W(t)$}
        |- (14,0) node (t-axis) [right] {$t$};
    \coordinate (a_1) at (2,0);
    \draw [dashed] (a_1) -- ++(0,3) coordinate (w_1);
    \draw (w_1) -- ++(-45:1.75) coordinate (a_2)
      edge[dotted] (intersection of w_1--a_2 and O--t-axis);
		\draw [dashed] (a_2) -- ++(0,3) coordinate (w_2);
		\draw (w_2) -- ++(-45:3.5) coordinate (a_3)
      edge[dotted] (intersection of w_2--a_3 and O--t-axis);
		\draw [dashed] (a_3) -- ++(0,3) coordinate (w_3);
		\draw (w_3) -- ++(-45:0.5) coordinate (a_4)
      edge[dotted] (intersection of w_3--a_4 and O--t-axis);
		\draw [dash dot] (a_4) -- ++(0,2) coordinate (w_4);
		\draw (w_4) -- ++(-45:3) coordinate (a_5)
      edge[dotted] (intersection of w_4--a_5 and O--t-axis);
		\draw[dash dot] (a_5) -- ++(0,1) coordinate (w_5);
    \coordinate (k) at (0,2.6);
    \draw[dotted] (k) node[left] {$k$}
        -- (t-axis |- k);
    \path (0,0) -- (k) node[midway, left] {$\mu_1$};
    \path (k) -- (w-axis) node[midway, left] {$\mu_2$};
    \node[draw=black, fill=white, rounded corners=2pt, above left=5] at (t-axis) {%
      \scriptsize
      \begin{tabular}{@{}r@{ }l@{}}
       \raisebox{2pt}{\tikz{\draw[dashed] (0,0) -- (5mm,0);}} & Exp($\mu_1$)\\
       \raisebox{2pt}{\tikz{\draw[dash dot] (0,0) -- (5mm,0);}} & Exp($\mu_2$)
      \end{tabular}
    };
\end{tikzpicture}
\caption{\label{fig:samplepath} Sample path of the VQT process $W(t)$ for the case $c=1$.}
\end{figure}
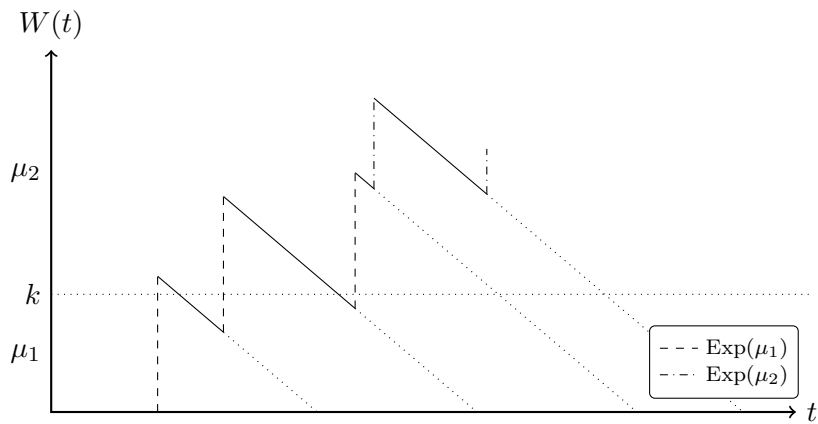

In case $c=1$, the process $\{ W(t), t \ge 0\}$ is sufficient for the state description. In fact, this process now corresponds to an M/M/1 queue in which the jump size depends on the state found upon arrival; see Figure~\ref{fig:samplepath} for an illustration of its sample path. In this case, the model is a special case of the M/G/1 variant of Model I in \cite{bekker2009levy}; see also e.g. \cite{gaver1962} for a classical related model with two service speeds.
Observe that for $c=1$ all matrices become scalars. In particular, using that $B_i = \tilde \Delta_i = \mu_i$ and $\tilde Q_i(x) = e^{-\mu_i x}$, and applying integration by parts, the integro-differential equations (\ref{eq:Q1}) and (\ref{eq:Q2}) in Theorem~\ref{th:dFx} can be written as 
\begin{eqnarray*}
F'(x) + \lambda \pi(0,0)(1-e^{-\mu_1 x})
 &=& \lambda \int_0^x e^{-\mu_1(x-y)} F'(y) dy + F'(0),\;\;\;  0 < x < k,  \\
F'(x) + \lambda \pi(0,0)(1-e^{-\mu_1 x})
 &=& \lambda \int_k^x e^{-\mu_2(x-y)} F'(y) dy + F'(0) \\
 & & + \lambda \int_0^k e^{-\mu_1(x-y)} F'(y) dy, \;\;\; x > k.  
\end{eqnarray*}
These equations can also be interpreted as level crossing equations, where the left and right hand sides correspond to the rate in and out of set $(0,x)$, respectively. 
Solving these equations, in terms of the density $F'(x)$, we obtain, for $0<x<k$, 
\[
F'(x) = \lambda \pi(0,0) e^{-(\mu_1-\lambda)x}, 
\]
whereas, for $x>k$, we have
\[ 
F'(x) = \lambda \pi(0,0) e^{-(\mu_1-\lambda)k}
\left[ \frac{(\mu_2-\mu_1)e^{-\mu_1 (x-k)}}{\mu_2-\mu_1-\lambda}  - \frac{\lambda e^{-(\mu_2-\lambda)(x-k)}}{\mu_2-\mu_1-\lambda} \right],
\]
with
\[
\pi(0,0) = \frac{(\mu_1/\lambda-1)(\mu_2/\lambda-1)}{(\mu_1/\lambda) (\mu_2/\lambda-1)-(\mu_2/\mu_1-1) e^{(\lambda -\mu_1) k}} \ . 
\]
The VQT density allows for an intuitive interpretation. Specifically, in the region $(0,k)$ jump sizes are always exp($\mu_1$), which implies that $F'(x)$ is proportional to the limiting workload density in an M/M/1 queue with service rate $\mu_1$ (and finite workload capacity $k$). Also, observe that sample paths of $W(t)$ in the region $(k,\infty)$ are always initiated by an upcrossing of $k$ with a jump of size exp($\mu_1$), after which all jumps are exp($\mu_2$) until a subsequent downcrossing of $k$. This implies that $W(t)$ in $(k,\infty)$ behaves as the workload process in an M/M/1 queue with service rate $\mu_2$, but with an exceptional first service time in a busy period that has rate $\mu_1$. This directly explains the mixture of the two exponentials in $(k,\infty)$.

\subsection{Numerical example for \texorpdfstring{$c=2$}{c=2}} \label{subs:c2}
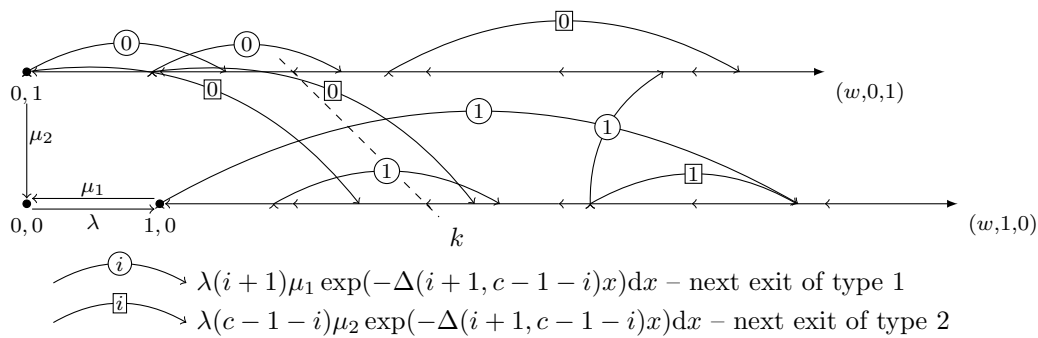
\begin{figure} 
    \centering
    \def\c{2}\input{diagram.tex}
    \caption{\label{fig:transitions} Sketch of the transition diagram for the case $c=2$.}
\end{figure}
\begin{figure}
    \centering
    \includegraphics[width=0.9\textwidth]{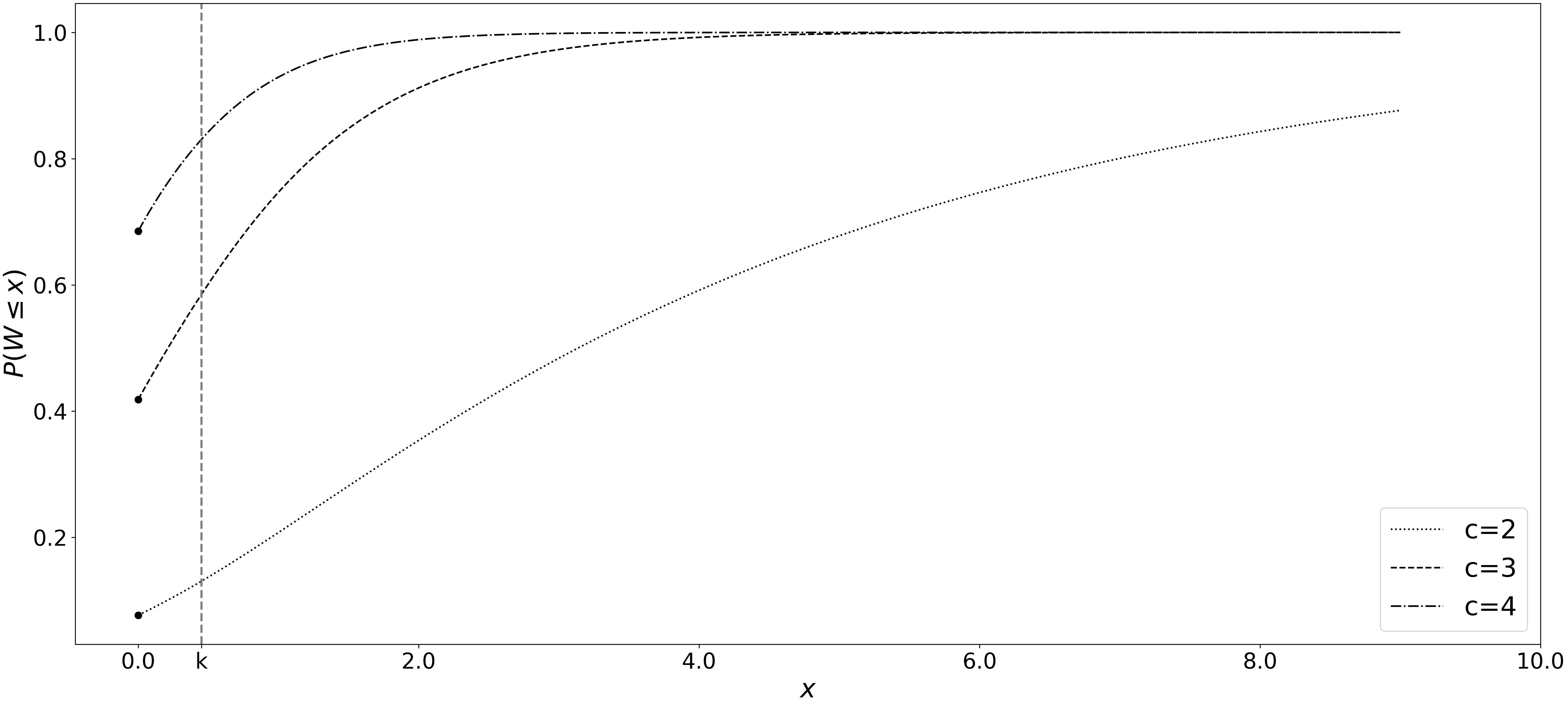}
    \caption[Plot of the Stationary VQT cumulative distribution]{\label{fig:serv2}Stationary VQT cumulative distribution function for $\lambda=2$, $\mu_1 = 0.75$, $\mu_2=1.12$ and $k=0.45$.\footnotemark[\fnsybm]}
\end{figure}

 
As an illustration of the balance equations and the limiting distribution, we consider the 2-server case in this section. 
A visual representation of the transition diagram of $(W(t), S_1(t), S_2(t))$ can be found in Figure~\ref{fig:transitions}.
To avoid excessive expressions, we focus on a numerical example with specific parameters.
\input{example_c2_l2_m0p75_M1p12_k0p45.tex}
We fix $k=\ex{k}$, $\lambda=\ex{lambda}$, $\mu_1=\ex{mu1}$ and $\mu_2=\ex{mu2}$. 
The second order differential equations of Theorem~\ref{th:dFx2} then looks as follows: 
for $0 \le x < \ex{k}$, we have
\begin{eqnarray*}
    F''_0(x) 
    \ex{A11[1,1]} F'_0(x) 
    \ex{A10[1,1]} F_0(x) 
    &=& 
    \ex{M01[1,1]} F'_0(0)
    \ex{M02[1,1]} \pi(0,1),
    \\
    F''_1(x) 
    \ex{A11[1,2]} F'_0(x) 
    \ex{A11[2,2]} F'_1(x) 
    \ex{A10[1,2]} F_0(x) 
    &=& 
    \ex{M01[1,2]} F'_0(0)
    \ex{M01[2,2]} F'_1(0)
    \\ & & 
    \ex{M02[1,2]} \pi(0,1)
    \ex{M02[2,2]} \pi(1,0),
\end{eqnarray*}
whereas, in the region  $x > \ex{k}$, we obtain
\begin{eqnarray*}
    \quad & & \hspace{-2cm}
    F''_0(x) 
    \ex{A21[1,1]} F'_0(x) 
    \ex{A21[2,1]} F'_1(x) 
    \ex{A20[2,1]} F_1(x) 
  = \\ & &
    \ex{M11[1,1]} F'_0(0)
    \ex{M11[2,1]} F'_1(0)
    \ex{M12[1,1]} \pi(0,1)
    \\ & &
    \ex{M13[1,1]} F_0(\ex{k}+)
    \ex{M13[2,1]} F_1(\ex{k}+) 
    \\ & &
    \ex{M21[1,1]}(F'_0(0)-F'_0(\ex{k}+))
    \\ & &
    \ex{M22[1,1]} \pi(0,1)
    \\ & &
    \ex{M23[1,1]} F_0(\ex{k}+)
    , \\ \quad & & \hspace{-2cm}
    F''_1(x) 
    \ex{A21[2,2]} F'_1(x) 
    \ex{A20[2,2]} F_1(x) 
  = \\ & & 
    \ex{M11[2,2]} F'_1(0)
    \ex{M12[2,2]} \pi(1,0)
    \ex{M13[2,2]} F_1(\ex{k}+)
    \\ & & +(\ex{M21[1,2]}) (F'_0(0)-F'_0(\ex{k}+))
    \\ & & \ex{M21[2,2]} (F'_1(0)-F'_1(\ex{k}+)) 
    \\ & & 
    \ex{M22[1,2]} \pi(0,1)
    \ex{M22[2,2]} \pi(1,0)
    \\ & & 
    +(\ex{M23[1,2]}) F_0(\ex{k}+).
\end{eqnarray*}
The boundary conditions in Corollary~\ref{cor:bc} are rather straightforward. 
By Theorem~\ref{th:FxMsol}, the solution of the above system of differential equations is unique depending on a constant $b_c$ and the components $\pi(1,0)$ and $\pi(0,1)$. 
That is the quantities $F'_0(0)$ and $F'_1(0)$ are determined given those values as follows
\begin{eqnarray*}
    F'_0(0) &=& 
    \ex{H16[1,1]} \pi(0,1)
    \ex{H16[2,1]} \pi(1,0) 
    \ex{psiH15[1]} b_c ,
    \hspace{1.35cm}  \\ 
    F'_1(0) &=& 
    \ex{H16[1,2]} \pi(0,1)
    \ex{H16[2,2]} \pi(1,0) 
    \ex{psiH15[2]} b_c .
    \hspace{1.25cm}  
\end{eqnarray*}
It follows that the remaining unknowns $\pi(0,0)$, $\pi(1,0)$, $\pi(0,1)$ and $b_c$
can be determined by imposing the following constraints
\begin{eqnarray*}
    F'_0(0) &=& 
    \ex{M3[1,1]} \pi(0,1),
    \\
    F'_1(0) &=& 
    \ex{M3[2,2]} \pi(1,0) 
    -\ex{lambda} \pi(0,0), \\
    0 &=& -\ex{lambda} \pi(0,0) + \ex{mu1} \pi(1,0) + \ex{mu2} \pi(0,1),\\
    1 &=& \pi(0,0) + \pi(1,0) + \pi(0,1) + F_0(\infty) + F_1(\infty),
\end{eqnarray*}
yielding
\begin{eqnarray*}
    \pi(0,0) &=& \ex{d0[1]} , \hspace{1.5cm}
    b_c = \ex{bc} ,\\
    \pi(0,1) &=& \ex{d1[1]} , \hspace{1.5cm}
    \pi(1,0) = \ex{d1[2]} .
\end{eqnarray*}
The expressions for $F_0(x)$ and $F_1(x)$, in the interval $0 \leq  x \leq \ex{k}$, are
\begin{eqnarray*}
    F_0(x) &=& 
    \ex{F1v1[1]} 
    \ex{F1v2[1]}
    \ex{F1v3[1]},\\
    F_1(x) &=& 
    \ex{F1v1[2]}
    \ex{F1v2[2]}
    \ex{F1v3[2]}, \\[0.4em]
& & \hspace*{-63pt}\textrm{whereas, in the interval $x > \ex{k}$, they are given by} \\[0.4em]
    F_0(x) &=& 
    \ex{F2v1[1]} 
    \\ & & 
    \ex{F2v2[1]}
    \ex{F2v3[1]},\\
    F_1(x) &=& 
    \ex{F2v1[2]}
    \ex{F2v2[2]}
    \ex{F2v3[2]}.
\end{eqnarray*}
Note that $\theta_1 = |\lambda - c\mu_1| = +0.5$ and $\beta_0 = \lambda - c\mu_2 = -0.24$.
The distribution of the virtual queueing time is visualized in Figure~\ref{fig:serv2}, along with the cases of 3 and 4 servers.%

\section{Numerical insights}
\label{sec:numerics}
In this section we focus on numerical insights. Specifically, we consider $W$, the stationary distribution of the queueing time or VQT. Note that the results in Section~\ref{sec:lim} contain more information, as they also provide the server state. To obtain VQT, observe that $P(W=0)=\sum_{i,j: 0\le i+j\le c-1} \pi(i,j)$ and $P(W \le x) = P(W=0) + F(x) \bone$.
From $F(x)$, we may also directly derive the mean stationary VQT, which we give here in matrix representation.
\footnotetext[\fnsybm]{\githubRemark}

\begin{lemma}\label{lm:mean.computation}
    The mean stationary VQT is computed as follows
    \begin{eqnarray}\label{eq:mean.computation}
        \mean[W] 
        &=& (F'(0) + \alpha_0 M_0 U^-_1) (U^+_1 - U^-_1)^{-1} I(0,k;U^+_1)
        \bone \nonumber \\ 
        & & - ((F'(0) + \alpha_0 M_0 U^-_1) (U^+_1 - U^-_1)^{-1} + \alpha_0 M_0) I(0,k;U^-_1)
        \bone \nonumber \\ 
        & & + (F(k) - b_c \psic - \alpha_1 M_1 - \alpha_2 (\tD_1-\tD_2) M_2) ((U^-_2)^{-1} - k I) \bone \nonumber \\
        & & - \alpha_2 (\tD_1^{-1} + k I) (\tD_1-\tD_2) M_2 \bone ,
    \end{eqnarray}
    where 
    $I(a,b;D)$ is defined as in \eqref{eq:int.x.exp(Dx)}.
\end{lemma}

\begin{proof} 
The results follow from the density $F'(x)$ and integration by parts, see \hyperref[pf:lemma.meanVQT]{Appendix}.
\end{proof}

\begin{figure}
\centering
\includegraphics[width=0.9\textwidth]{{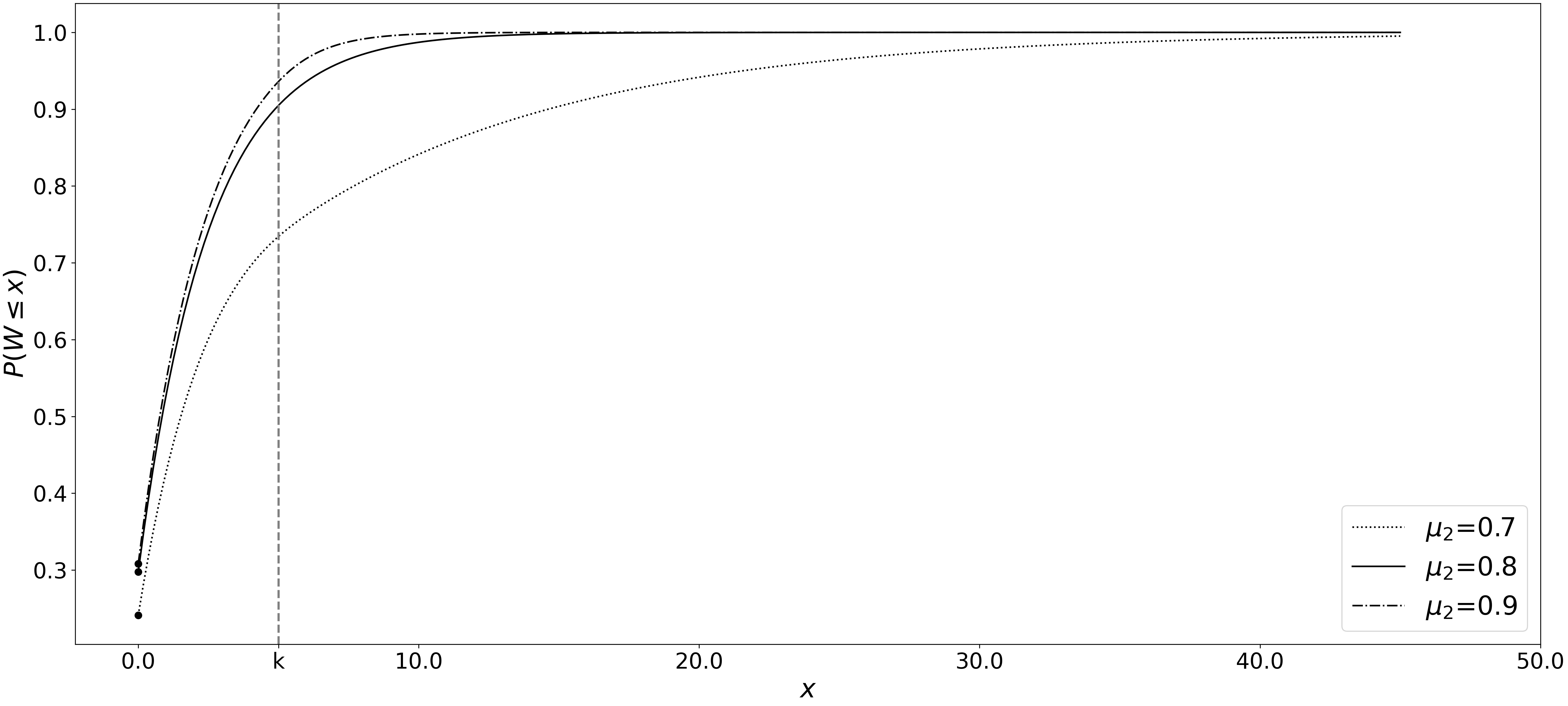}}
\caption[Plot of the Stationary VQT cumulative distribution]{\label{fig:RENE.A2} Cdf of the stationary VQT for $k=5$, $c=3$, $\lambda = 2$, $\mu_1 =0.8$, and $\mu_2\in\{0.7,0.8,0.9\}$.\footnotemark[\fnsybm]}
\end{figure}

First, we consider the impact of having different $\mu_1$ and $\mu_2$. 
In case $\mu_1 = \mu_2$, the system corresponds to the classical M/M/c queue, whereas there is a slowdown (speedup) effect when $\mu_2 < \mu_1$ ($\mu_2 > \mu_1$). 
As a basic example, we take $k=5$, $c=3$, $\lambda=2$, and $\mu_1=0.8$, whereas $\mu_2 \in \{0.7, 0.8, 0.9\}$.
The cumulative distribution function (cdf) of $W$ is visualized in Figure~\ref{fig:RENE.A2}. Clearly, in case of a slowdown ($\mu_2 = 0.7$), the queueing time strongly deteriorates compared to the standard situation where $\mu_2=\mu_1=0.8$. In fact, if $\mu_2 \le 2/3$ the system would even become unstable. For the current example, the impact of a speedup ($\mu_2 = 0.9$) is relatively small compared to the standard situation, as the basic service rate of 0.8 is already sufficient to provide reasonable queueing times. 
Hence, taking differences in service rates into account is crucial to provide reliable queueing times, especially in case of slowdowns.

\begin{figure}
\centering
\includegraphics[width=0.9\textwidth]{{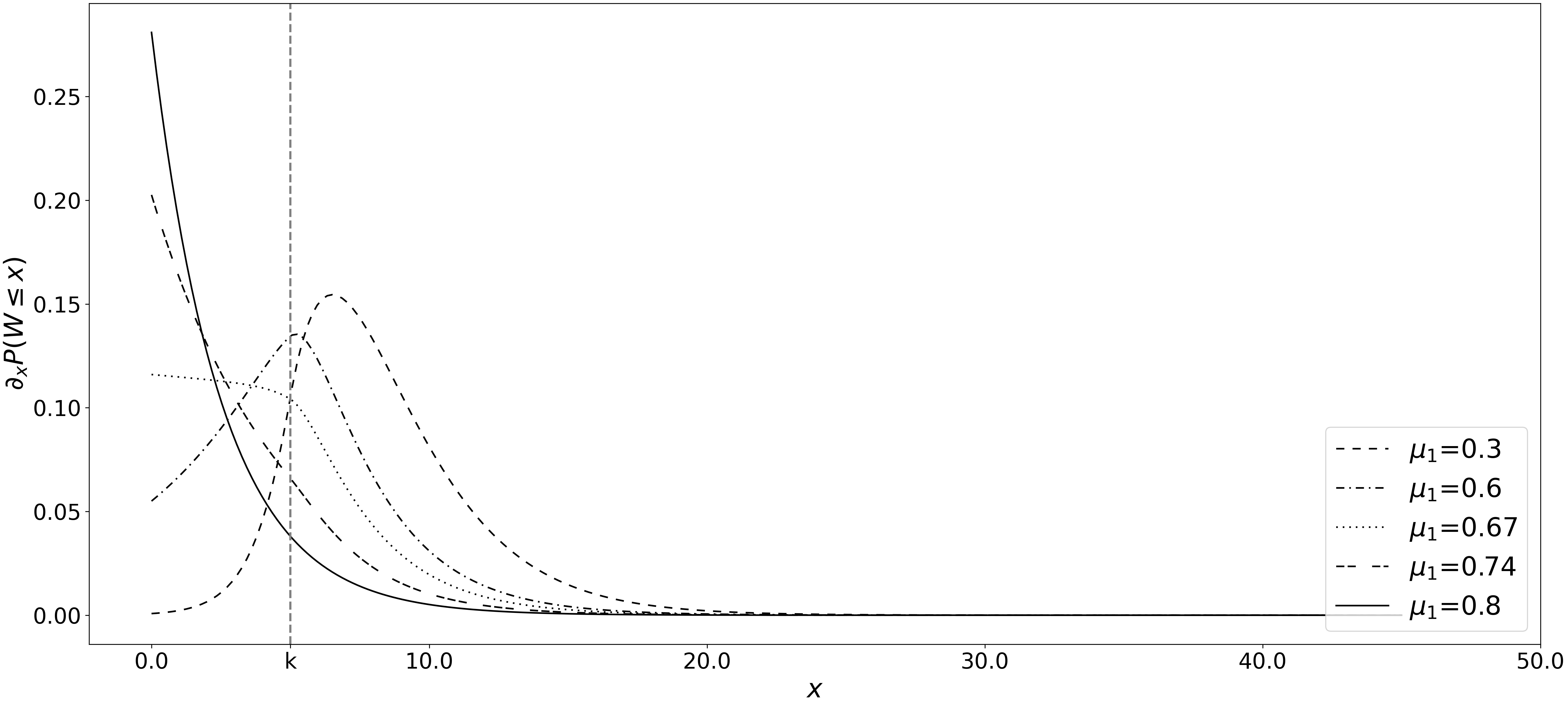}}
\caption[Plot of the Stationary VQT density function]{\label{fig:RENE.C2}Stationary VQT density for $k=5$, $c=3$, $\lambda = 2$, $\mu_2 =0.8$, and $\mu_1\in\{0.3,0.6,0.67,0.74,0.8\}$.\footnotemark[\fnsybm]}
\end{figure}

Second, the shape of the VQT density may also be strongly affected by speedups (or slowdowns), i.e., differences in $\mu_1$ and $\mu_2$. The VQT density is strictly decreasing for the standard M/M/c queue, which will also hold in case $\mu_2<\mu_1$ (slowdown).
However, this is no longer necessarily the case for speedups, see Figure~\ref{fig:RENE.C2} for $k=5$, $c=3$, $\lambda=2$, $\mu_2=0.8$, and $\mu_1 \in \{0.3, 0.6, 0.67, 0.74, 0.8\}$. 
In particular, for more extreme variants of a speedup effect, the peak in the VQT density may be around or above level $k$. 

\begin{figure}
\centering
\includegraphics[width=0.9\textwidth]{{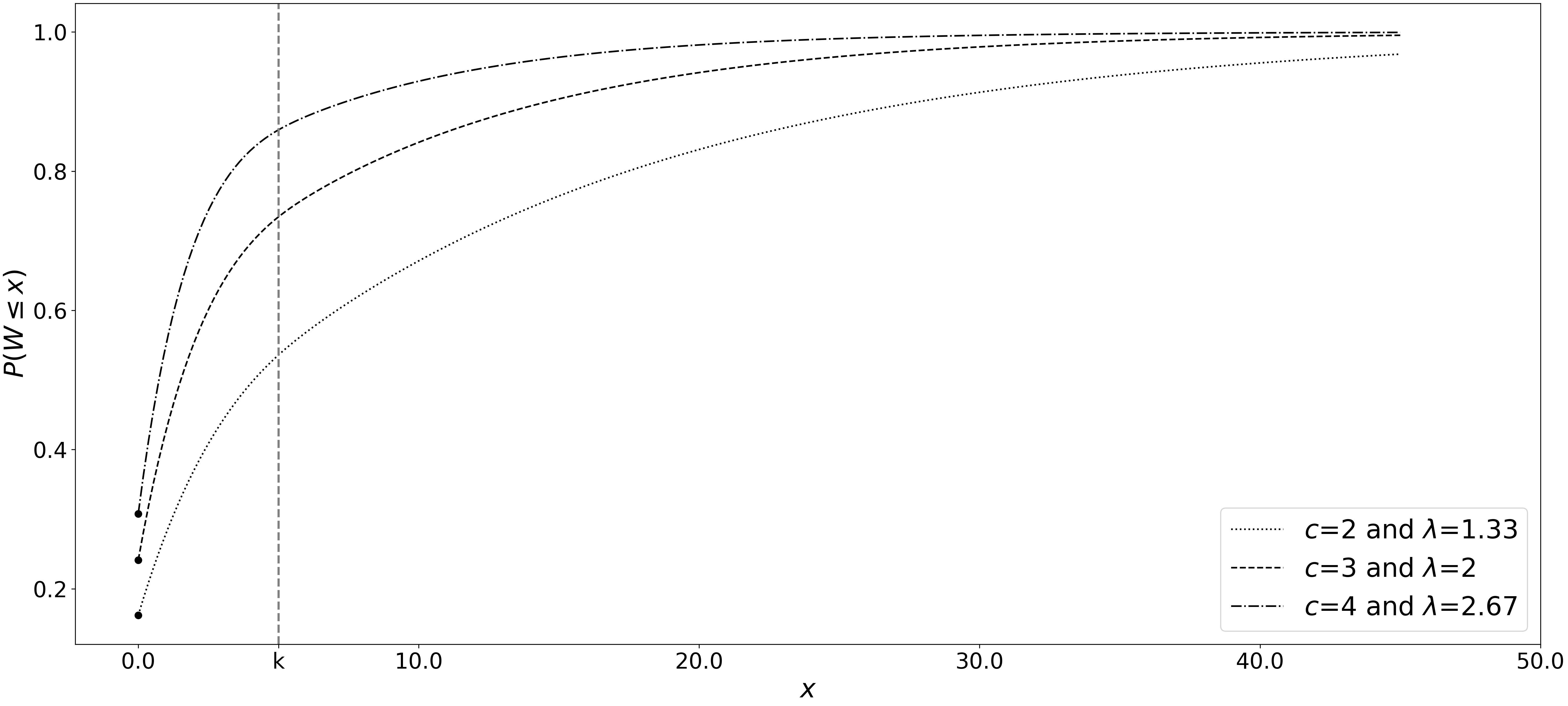}}
\caption[Plot of the Stationary VQT cumulative distribution]{\label{fig:RENE.B2} Cdf of the stationary VQT for $k=5$, $\mu_1 =0.8$, $\mu_2 =0.7$, $(c,\lambda)\in\{(2,4/3),(3,2),(4,8/3)\}$.\footnotemark[\fnsybm]}
\end{figure}

Third, we consider the impact of the number of servers, i.e., scale of the system. Let $k=5$, $\mu_1=0.8$, $\mu_2=0.7$ (slowdown), and consider systems with 2, 3, and 4 servers. We let $\lambda$ be 4/3, 2, and 8/3, respectively, such that the loads $\lambda/(c\mu_i)$, for $i=1,2$ are identical. The cdf of the stationary VQT is presented in Figure~\ref{fig:RENE.B2}. 
Clearly, as the number of servers increases the queueing time improves, which is in line with economies of scale for regular M/M/c queues. We like to note that the relative ordering of cdf's below $k$ may change in case $\lambda \ge c \mu_1$ (which we did not visualize here).

\begin{figure}
    \centering
    \includegraphics[width=0.9\textwidth]{{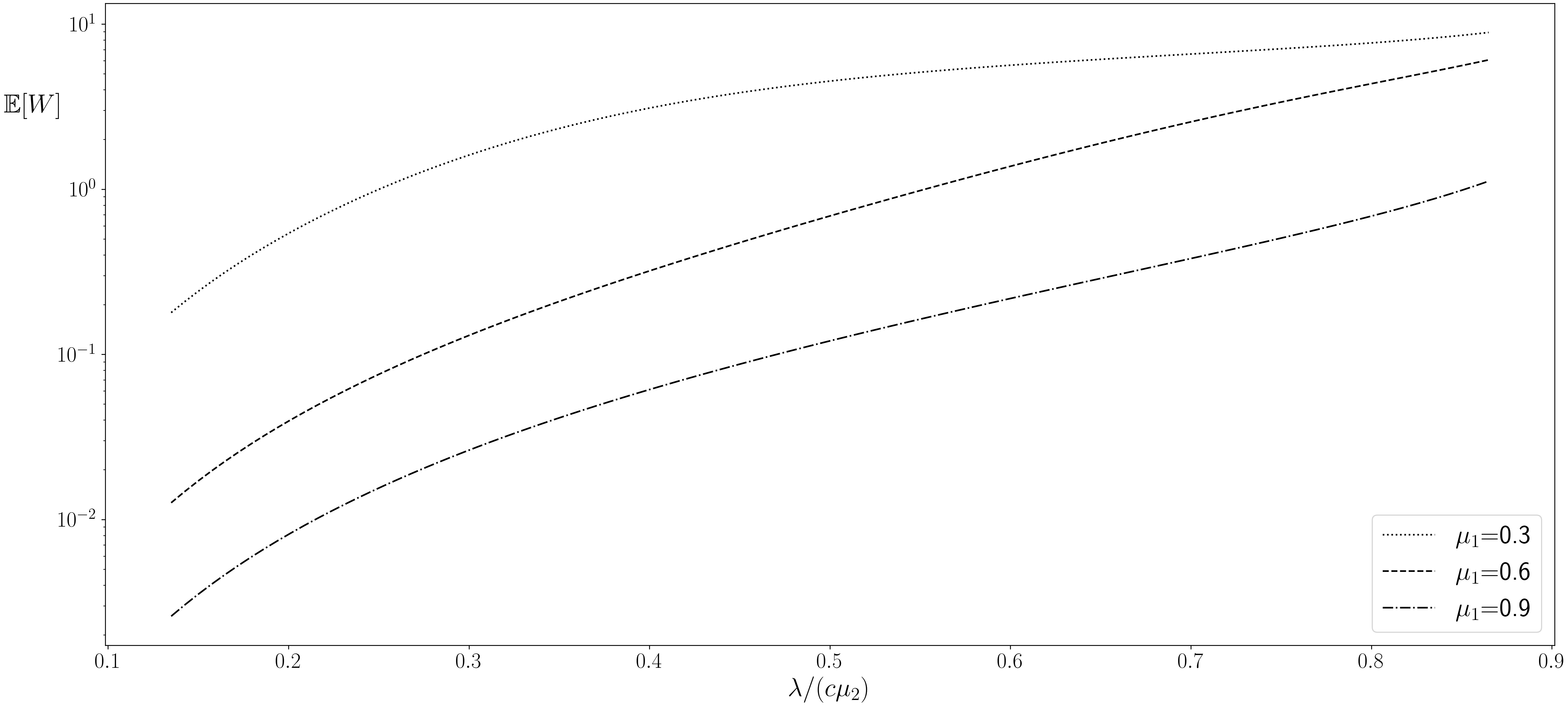}}
    \caption[Plot of the Stationary Expected VQT]{\label{fig:RENE.D2}Expected VQT for $k=5$, $\mu_2 = 0.8$, $c = 3$ and $\mu_1\in\{0.3,0.6,0.9\}$.\footnotemark[\fnsybm]}
\end{figure}

Finally, we consider the expected VQT. 
It is well known that $\mean[W]$ is convex and increasing in $\lambda$ for regular M/M/c queues with $c$ and $\mu$ fixed. This property is not necessarily preserved in the current model, see Figure~\ref{fig:RENE.D2}. Specifically, in case $\mu_1$ is relatively small ($\mu_1=0.3$ in Fig.~\ref{fig:RENE.D2}), a large fraction of the customers will experience a VQT of around $k$, assuming the system to be stable. This destroys the convexity of $\mean[W]$ as a function of $\lambda$. Moreover, the impact of $\mu_1$ is also considerable for more heavily loaded systems. For instance, comparing $\mean[W]$ for different $\mu_1 \in \{0.3,0.6,0.9\}$ with fixed $\lambda/(c\mu_2)=0.9$, we see that the mean VQT is much smaller when a lot of customers can be served with rate $\mu_1$ (i.e., for $\mu_1=0.9$).
To conclude this section, we note that neglecting the differences in service rate leads to rather inadequate performance characteristics.

\appendix

\section{Proofs} \label{app:proofs}
In this appendix, we present the technical proofs of the results presented throughout the paper.

\begin{proofof}{{Theorem~\ref{th:dFx}}}\label{pf:th:dFx}
Let $0 \le x < k$, and $0 \le i \le c-1$ be fixed. 
Define $P_i(W(t)\in A) = P(W(t) \in A, S(t)=(i,c-1-i))$, with $A\subset\bR$. Conditioning on the jump size being exactly equal to $x-y$, 
we get 
\begin{eqnarray} \label{F_bef_k}
    F_i(t,x)
    &=&  \int_{0}^x P_i(W(t-h) \leq y+h) \lambda h (i+1)\mu_1 e^{-\Delta(i+1,c-1-i)(x-y)} dy \nonumber\\
    & & + \int_{0}^x P_{i-1}(W(t-h) \leq y+h) \lambda h (c-i)\mu_2 e^{-\Delta(i,c-i)(x-y)} dy \nonumber\\
    & & + (1-\lambda h) P_i(h < W(t-h) \le x+h) + o(h) \ . 
\end{eqnarray}

Define
\[ \mQ_\kappa(x) = \exp(-(\mu_\kappa I+ \mD_{c-1}) x) , \;\;  \kappa \in \{1, 2\}, \]
let $t \rightarrow \infty$, divide by $h$, and rearrange terms to get
\begin{eqnarray*}
    \lefteqn{\frac{1}{h} (F_i(x)-F_i(x+h)) =  -\lambda F_i(x+h) - (1-\lambda h) \frac{1}{h}(F_i(h)-0) } \\
    & &  + \lambda \sum_{j=0}^{c-1} \int_{0}^x (\pi(j,c-1-j) + F_j(y+h)) \sum_{k=0}^{c-1} \mQ_1(x-y) B_1(k,i) dy \\
    & &  + o(h)/h \ .
\end{eqnarray*}
Letting $h \rightarrow 0$ and multiplying both sides  by $-1$, and noting that $F(0)=0$, 
we get Equation~\eqref{eq:Q1}, after noticing that $\mQ_1 B_1 = B_1 \tQ_1$.\\

For $x > k$, again with $P_i(W(t)\in A) = P(W(t) \in A, S(t)=(i,c-1-i))$ for $A\subset\bR$, we have
\begin{eqnarray}  \label{F_aft_k}
\lefteqn{F_i(t,x) = } \nonumber\\
&  & \int_{0}^k P_i(W(t-h) \le y+h) \lambda h (i+1)\mu_1 e^{-\Delta(i+1,c-1-i)(x-y)} dy \nonumber\\
& &+ \int_{0}^k P_{i-1}(W(t-h) \le y+h) \lambda h (c-i)\mu_2 e^{-\Delta(i,c-i)(x-y)}dy \nonumber\\
& & + \int_{k}^x P_i(W(t-h) \le k) \lambda h (i+1)\mu_1 e^{-\Delta(i+1,c-1-i)(x-y)}dy \nonumber\\
& & +\int_{k}^x P_{i-1}(W(t-h) \le k) \lambda h (c-i)\mu_2 e^{-\Delta(i,c-i)(x-y)} dy \nonumber\\
& & +\int_{k}^x P_i(k < W(t-h) \le y+h) \lambda h (c-i)\mu_2 e^{-\Delta(i,c-i)(x-y)} dy \nonumber\\
& & + \int_{k}^x P_{i+1}(k < W(t-h) \le y+h) \lambda h (i+1)\mu_1 e^{-\Delta(i+1,c-1-i)(x-y)}dy \nonumber\\
& & + (1-\lambda h) P_i(h<W(t-h) \le x+h) 
    + o(h) \ .
\end{eqnarray}
Following the same steps as above, we get
\begin{eqnarray*}
    \lefteqn{-F'(x) = - \lambda F(x) - F'(0)} \\
& & + \lambda \int_0^k (\state{c-1}+ F(y)) \mQ_1(x-y) B_1 dy + \lambda (\state{c-1}+F(k)) \int_k^x \mQ_1(x-y) B_1 dy \\
&  & + \lambda \int_k^x (\state{c-1}+ F(y)) \mQ_2(x-y) B_2 dy
- \lambda (\state{c-1}+F(k)) \int_k^x \mQ_2(x-y) B_2 dy \ .
\end{eqnarray*}
Substituting $\mQ_\kappa B_\kappa = B_\kappa \tQ_\kappa$, $\kappa\in\{1,2\}$, yields Equation~\eqref{eq:Q2}. 
\end{proofof}

\begin{proofof}{Theorem~\ref{th:dFx2}}\label{pf:th:dFx2}
     First consider Equation~\eqref{eq:Q1}. Note that 
\[ \tQ_\kappa'(x) = - \tQ_\kappa(x) \tD_\kappa, \;\;\; \kappa = 1,2.\]
Since $B_1$ is invertible, we can use Equation~\eqref{eq:Q1} to get, for $0 \le x \le k$,
\begin{eqnarray}
    \lambda \int_0^x F(y) B_1 \tQ_1(x-y) dy
    &=& - F'(x) + F'(0) + \lambda F(x) \nonumber \\
    & & - \lambda \state{c-1} B_1 (I - \tQ_1(x)) \tD_1^{-1}. \label{eq:lint}
\end{eqnarray}
Taking the derivative with respect to $x$ on both sides of Equation~\eqref{eq:Q1}, we get 
\[ 
F''(x) =  
 \lambda F'(x)  
- \lambda F(x) B_1 
+ \lambda \int_0^x F(y) B_1 \tQ_1(x-y) \tD_1 dy
-\lambda \state{c-1} B_1 \tQ_1(x).
\]
Substituting Equation~\eqref{eq:lint} in the above equation we get Equation~\eqref{eq:F12} with
\begin{eqnarray}
    \alpha_0 &=&  F'(0) \tD_1 - \state{c-1} \lambda B_1, \nonumber\\
        &=&   [\state{c-1}(\lambda I + \mD_{c-1}) - \state{c-2} \lambda \hat I]\tD_1 - \state{c-1} \lambda B_1 \ . 
\end{eqnarray}
Here we have used Equation~\eqref{F.zero.cond} to eliminate $F'(0)$. This is a second order linear differential equation with constant coefficients and a constant driving function on the right hand side.

Next we consider Equation~\eqref{eq:Q2}. First, note that, for $x > k$, we get, by applying Equation~\eqref{eq:lint}, that
\begin{eqnarray}
    \lambda \int_0^k F(y) B_1 \tQ_1(x-y) dy
    &=& \left[\lambda \int_0^k F(y) B_1 \tQ_1(k-y) dy \right] \tQ_1(x-k) \nonumber \\
    & = & [- F'(k) + F'(0) + \lambda F(k) \nonumber \\
    & & 
    - \lambda \state{c-1} B_1 (I - \tQ_1(k)) \tD_1^{-1}] \tQ_1(x-k) .\label{eq:iok}
\end{eqnarray}
We also have, for $x > k$,  
\begin{equation} \label{eq:ikx}
\int_k^x B_\kappa \tQ_\kappa(x-y) dy = B_\kappa (I - \tQ_\kappa(x-k)) \tD_\kappa^{-1}  , \;\; \kappa=1,2.
\end{equation}
Substituting in the RHS of Equation~\eqref{eq:Q2} we get, for $x > k$, 
\begin{eqnarray}
F'(x) & = & \lambda F(x) - \lambda \int_k^x F(y) B_2 \tQ_2(x-y) dy \nonumber \\
& & + \alpha_1 \tD_2^{-1} -\alpha_2 \tQ_1(x-k) - \lambda F(k) B_2 \tQ_2(x-k) \tD_2^{-1},\label{eq:vfk}
\end{eqnarray}
where 
\begin{eqnarray*}
\alpha_1 \tD_2^{-1} & = & F'(0) - \lambda F(k) (B_1\tD_1^{-1} - B_2\tD_2^{-1}) - \state{c-1} \lambda B_1\tD_1^{-1}, \\
\alpha_2 & = & \lambda F(k)(I - B_1\tD_1^{-1}) - (F'(k)-F'(0)) - \state{c-1} \lambda B_1 \tD_1^{-1}.
\end{eqnarray*}
Differentiating both sides of Equation~\eqref{eq:vfk}, we get 
\begin{eqnarray}
F''(x) & = & \lambda F'(x) - \lambda F(x)B_2 +\lambda \int_k^x F(y) B_2 \tQ_2(x-y)\tD_2 dy \nonumber \\
& &  \alpha_2 \tQ_1(x-k)\tD_1 + \lambda F(k) B_2 \tQ_2(x-k).\label{eq:vfk2}
\end{eqnarray}
Using \eqref{eq:vfk} we have 
\begin{eqnarray}
    \lambda \int_k^x F(y) B_2 \tQ_2(x-y) dy 
    & = & - F'(x) + \lambda F(x) + \alpha_1 \tD_2^{-1}  - \alpha_2 \tQ_1(x-k)  \nonumber\\
    & & - \lambda F(k) B_2 \tQ_2(x-k) \tD_2^{-1}. \label{eq:vfk3}
\end{eqnarray}
Substituting  Equation~\eqref{eq:vfk3} in the RHS of  Equation~\eqref{eq:vfk2} we get Equation~\eqref{eq:F22}. 
This completes the proof.
\end{proofof}

\begin{proofof}{Corollary~\ref{cor:bc}}\label{pf:cor:bc}
Equation~\eqref{Fzero} follows from the definition of $F$. Equation~\eqref{F.cont.cond} follows by taking the left and right limits at $k$ in Equations \eqref{F_bef_k} and \eqref{F_aft_k}, respectively.
Equation~\eqref{F.diff.cond} follows by taking the left and right limits at $k$ in Equations \eqref{eq:Q1} and \eqref{eq:Q2}, respectively.

The balance equation for state $(0,i,c-1-i)$ yields
\[ (\lambda + i\mu_1 + (c-1-i)\mu_2) \pi(i,c-1-i) = F_i'(0) + \lambda \pi(i-1,c-1-i), \;\;\; 0 \le i \le c-1.\]
In matrix form this can be written as 
\[
    F'(0) = \state{c-1}(\lambda I + \mD_{c-1}) - \state{c-2} \lambda \hat I ,
\]
 which is Equation~\eqref{F.zero.cond}.
\end{proofof}

\begin{lemma}\label{lem:M.matrix} Let $M$ be a $c \times c$ matrix with entries given by
    \begin{eqnarray*}
        M(0,c-1) & = & c \mu_2, \\
        M(1,c-1) & = & \phantom{c} \mu_1, \\
        M(i,i-1) & = & \frac{i}{c-i} \frac{\mu_1}{\mu_2}, \;\;\; 0 < i < c-1,\\
        M(i,j)& = & 0 \;\;\;\mbox{for all other } (i,j) .
    \end{eqnarray*}
    It is non-singular and satisfies the following equation
   \begin{equation}\label{eq:M} 
        M (B_1 - \mu_1 I  - \mD_0) = (B_2 - \mu_2 I  - \mD_0) \ .
   \end{equation}
\end{lemma}
\begin{proof} For $0 \leq i,j\leq c-1$, we write 
    \begin{eqnarray*}
          B_1(i,j) & = & \Kronecker{i}{j} (i + 1) \mu_1 + \Kronecker{i+1}{j} (c - i - 1) \mu_2, \\
          B_2(i,j) & = & \Kronecker{i}{j} (c - i) \mu_2 + \Kronecker{i-1}{j} i \mu_1, \\
        \mD_0(i,j) & = & \Kronecker{i}{j} (i \mu_1 + (c - i - 1) \mu_2), \\
            M(i,j) & = & \Kronecker{i + c - 1}{j} c \mu_2 + \Kronecker{i + c - 1}{j + 1} i \mu_1 \\ 
                   &   & - \Kronecker{i}{j + 1} i \mu_1 / ((c - i)\mu_2 ),
    \end{eqnarray*}
    where $\Kronecker{i}{j}$ is the Kronecker delta function, that is 
    \begin{equation}
        \Kronecker{i}{j} \defeq \left\{\begin{array}{ll}
            1 & \mbox{if } i=j , \\
            0 & \mbox{otherwise.}
        \end{array}\right.
    \end{equation}
    Let $Y_\kappa=B_\kappa-\mu_\kappa I -\mD_0$ , $\kappa\in\{1,2\}$. It follows that 
    \begin{eqnarray*}
     Y_1(i,j) & = & (\Kronecker{i+1}{j} - \Kronecker{i}{j}) (c - i - 1) \mu_2 , \\
     Y_2(i,j) & = & (\Kronecker{i-1}{j} - \Kronecker{i}{j}) i \mu_1 . 
    \end{eqnarray*}
    Then by matrix multiplication we have
    \begin{eqnarray*}
        (M Y_1)(i,j) 
        & = & \sum_{k=0}^{c-1} M(i,k) Y_1(k,j) \\
        & = & \sum_{k=0}^{c-1} \Big( \Kronecker{i + c - 1}{k} c \mu_2 + \Kronecker{i + c - 1}{k + 1} i \mu_1 \Big) Y_1(k,j) \\ 
        & & - \sum_{k=0}^{c-1} \Kronecker{i}{k + 1} \frac{i}{c-i} \frac{\mu_1}{\mu_2} Y_1(k,j) \\
        & = & \Big( \Kronecker{i}{0} c \mu_2 + \Kronecker{i}{1} \mu_1 \Big) Y_1(c-1,j)
              - \frac{i}{c-i} \frac{\mu_1}{\mu_2} Y_1(i-1,j) \\
        & = & - \frac{i}{c-i} \frac{\mu_1}{\mu_2} Y_1(i-1,j) 
        = (\Kronecker{i-1}{j} - \Kronecker{i}{j}) i \mu_1 = Y_2(i,j) .
    \end{eqnarray*}
\end{proof}

\begin{proofof}{Lemma~\ref{lm:right.eigenvec}}\label{pf:lm:right.eigenvec}
    Let us assume that the left equation in \eqref{eq:right.eigenvec} holds, that is 
    $\alpha_1 \cdot \tpsic = 0$. 
    Then using the definition of $\alpha_1$, given in \eqref{eq:alp1}, 
    we have that 
    \[ \alpha_0 \tD_1^{-1} \tD_2  \tpsic - \lambda F(k+) (B_1\tD_1^{-1}\tD_2  - B_2)  \tpsic = 0 \]
    and we are going to show that 
    \begin{equation}\label{eq:tpsic.to.show}
    (B_1 - \tD_1)  \tD_1^{-1} \tD_2  \tpsic = 0 ,
    \end{equation}
    so that the column vector $\tD_1^{-1} \tD_2  \tpsic$ 
    is parallel to $\tphic$, 
    because it is a right eigenvector of the matrix $(B_1 - \tD_1)$ 
    corresponding to the null eigenvalue, whose multiplicity is one.

    If Equation \eqref{eq:tpsic.to.show} holds, we have that
    \begin{eqnarray}
        (B_1 - \tD_1)  \tD_1^{-1} \tD_2  \tpsic  
        & = & (B_1 \tD_1^{-1} \tD_2 - \tD_2)  \tpsic \nonumber\\
        & = & (B_1 \tD_1^{-1} \tD_2 - B_2)  \tpsic = 0, \label{eq:tpsic.implication}
    \end{eqnarray}
    where in the last equality we used the relation $\tD_2 \tpsic = B_2 \tpsic$ given by \eqref{eq:tpsic.prop}.

    Equation \eqref{eq:tpsic.implication}, together with \eqref{eq:tpsic.to.show}, 
    implies that $\alpha_0 \tD_1^{-1} \tD_2  \tpsic = 0$ and therefore it also implies the result. 

    To prove \eqref{eq:tpsic.to.show}, we continue from \eqref{eq:tpsic.implication} by rewriting it in the following way
    \begin{eqnarray*}
        (B_1 \tD_1^{-1} \tD_2 - B_2) \tpsic 
        & = & (B_1 \tD_1^{-1} B_2 - B_2) \tpsic \\
        & = & (B_1 \tD_1^{-1} - I) B_2 \tpsic \\
        & = & ((\mu_1 I + \mD_{c-1})^{-1} B_1 - I) B_2 \tpsic \\
        & = & (\mu_1 I + \mD_{c-1})^{-1} ( B_1 - \mu_1 I - \mD_{c-1}) B_2 \tpsic \\
        & = & (\mu_1 I + \mD_{c-1})^{-1} M^{-1} M ( B_1 - \mu_1 I - \mD_{c-1}) B_2 \tpsic \\
        & = & (\mu_1 I + \mD_{c-1})^{-1} M^{-1} (B_2 - \mu_2 I  - \mD_{c-1}) B_2 \tpsic \\
        & = & (\mu_1 I + \mD_{c-1})^{-1} M^{-1} B_2^{-1} (B_2 - \tD_2) \tpsic = 0.
    \end{eqnarray*}  
    In the first equality we used \eqref{eq:tpsic.prop}, 
    in the third one that $B_1 \tD_1^{-1} = (\mu_1 I + \mD_{c-1})^{-1} B_1$, given by the definition \eqref{def:tD.kappa},
    in the sixth one we used the result of Lemma \ref{lem:M.matrix}, where the matrix $M$ is defined.
    Finally in the last two equalities we use again the definition \eqref{def:tD.kappa} and the hypothesis \eqref{eq:tpsic.prop}.
\end{proofof}

\begin{proofof}{Theorem~\ref{th:Fx2sol}}\label{pf:th:Fx2sol}
We consider the two regions of $x$ separately. First assume that $0 < x < k$. 
Equation~\eqref{eq:F12} is a non-homogeneous linear system of  ordinary second order differential equations. Hence, we first try a homogeneous solution of the type
\[ F_h(x) = e^{\theta x}\phi,\]
where $\phi$ is a row vector of length $c$. Substituting in $\cL_1 F_h(x) = 0$
and cancelling $e^{\theta x}$, we get Equation~\eqref{eq:eig1}, with $2c$ solutions $\{(\theta_i,\phi_i), 0 \le i \le 2c-1\}$. 
The eigenvalues $\theta_i$'s are given in Equations~\eqref{eq.theta.neg} and \eqref{eq.theta.pos}.
The homogeneous solution to Equation~\eqref{eq:F12} is then given by 
\begin{equation} \label{eq:solh1}
F_h(x) = \sum_{i=0}^{2c-1} a_i e^{\theta_i x} \phi_i, \;\;\; 0 \le  x \le k,
\end{equation}
where the $2c$ constants $\{a_i, 0 \le i \le 2c-1\}$ are to be determined.

For the particular solution we should look for a function of the following type 
\begin{equation}\label{eq:part.sol}
    F_p(x) = \eta x + \zeta,
\end{equation}
because $\theta_{c-1}=0$ is an eigenvalue for the homogeneous solution.
By substitution in \eqref{eq:F12} we get that the following equation has to be satisfied 
\begin{equation}\label{eq:F12.algebraic} 
    -\eta (\lambda I - \tD_1) 
    + \eta x \lambda (B_1 - \tD_1)
    + \zeta \lambda (B_1 - \tD_1) = \alpha_0 ,
\end{equation}
where the vector $\zeta$ can be chosen such that $\zeta \cdot \phic =0$, because for all $a \in \bR$, 
\[(\zeta + a \, \phic ) \lambda (B_1 - \tD_1) = \zeta \lambda (B_1 - \tD_1). \]
In addition, in order to have \eqref{eq:F12.algebraic} satisfied for any $x$, the coefficient of the linear term should be null implying that $\eta = a \, \phic$.

Taking the scalar product of both sides of \eqref{eq:F12.algebraic} by the right eigenvector $\tphic$ satisfying \eqref{eq:tphic.prop}, we get that
\[-a \, \phic (\lambda I - \tD_1)  \cdot \tphic = \alpha_0 \cdot \tphic , \]
implying that 
\begin{equation}\label{eq:a}
    a = -\frac{\alpha_0 \cdot \tphic}{\phic (\lambda I - \tD_1)  \cdot \tphic} . 
\end{equation}
Note that the linear term is missing if $\alpha_0 \cdot \tphic=0$.

To derive the value of $\zeta$ we rewrite the equation $\zeta \cdot \phic=0$ in matrix form as $\zeta \, \diag(\phic) = 0$.
By adding this equation to \eqref{eq:F12.algebraic} we get 
\[ -\eta (\lambda I - \tD_1)  + \zeta (\lambda (B_1 - \tD_1) + \diag(\phic)) = \alpha_0 \]
and since $(\lambda (B_1 - \tD_1) + \diag(\phic))$ is not singular, we have that 
\[ \zeta = \alpha_0 M_0+ a \, \phic  (\lambda I - \tD_1) M_0 ,\]
where $M_0$ is as given in Equation~\eqref{eq:M0}.

Finally, we have that the particular solution is equal to 
\begin{equation}\label{eq:part.sol.[0,k]}
    F_p(x) = a \, \phic \left( x I + (\lambda I - \tD_1) M_0 \right)+ \alpha_0 M_0 , \;\;\; 0 \le x \le k ,
\end{equation}
where $a$ is defined as in Equation \eqref{eq:a}.
Below we show that $\alpha_0 \cdot \tphic=0$, implying that $a=0$ and thereby Equation \eqref{eq:gsol1}.

Next consider the region $x > k$, where $F$ satisfies Equation~\eqref{eq:F22}. 
It is also a non-homogeneous linear system of  ordinary second order differential equations. 
As before we try a homogeneous solution of the type
\[ F_h(x) = e^{\beta x}\psi,\]
where $\psi$ is a row vector of length $c$. 
Substituting in $\cL_2 F_h(x) = 0$ and cancelling $e^{\beta x}$ we get Equation~\eqref{eq:beta}, which has  $2c$ solutions $(\beta_i,\psi_i), \; 0 \le i \le 2c-1$. The eigenvalues $\beta$'s are given in Equations~\eqref{eq.beta.neg} and \eqref{eq.beta.pos}. 
The homogeneous solution to Equation~\eqref{eq:F22} is then given by 
\begin{equation} \label{eq:solh2.guess}
F_h(x) = \sum_{i=0}^{2c-1} b_i e^{\beta_i (x-k)} \psi_i, \;\;\; x \ge k.
\end{equation}
Since the solution has to be bounded we immediately get that 
$b_i = 0$ for $ c < i \le 2c-1$, since the corresponding $\beta_i$'s are strictly positive. Moreover, we have $\beta_c=0$ and $\psic$ is as given in Equation~\eqref{eq:psic}.
It follows that the homogeneous solution to Equation~\eqref{eq:F22} can be written as
\begin{equation} \label{eq:solh2}
    F_h(x) = \sum_{i=0}^{c-1} b_i e^{\beta_i (x-k)} \psi_i + b_c \psic, \;\;\; x \ge k,
\end{equation}
where the $c+1$ constants $\{b_i, 0 \le i \le c\}$ are to be determined.

Next we determine the particular solution.
Similarly to what we have done before for the interval $[0,k]$, 
    the particular solution associated with the constant term  $\alpha_1$ in the right hand side of \eqref{eq:F22} 
    would be of the form \eqref{eq:part.sol}, since the associated homogeneous equation $\cL_2 F_h(x) = 0$
    admits the constant function as solution.
    However, in this case, the boundary condition, requiring $\lim_{x\to\infty} F(x)$ to be bounded,
    implies that the vector $\eta$ is zero and therefore that $\alpha_1 \cdot \tpsic = 0$.
    By applying Lemma \ref{lm:right.eigenvec}, this also implies that the linear term in \eqref{eq:part.sol.[0,k]} is missing.

Eventually it follows that the particular solution in the region $x>k$ is given by
\begin{equation}\label{eq:solp2}
    F_p(x) = \alpha_1 M_1 + \alpha_2 \tQ_1(x-k) (\tD_1-\tD_2) M_2  , \;\;\; x \ge k,
\end{equation}
as can be verified by direct substitution,
where  $M_1$ and $M_2$ are as given in Equations~\eqref{eq:M1} and \eqref{eq:M2}.
The general solution is then as given in Equation~\eqref{eq:gsol2}. This completes the proof.
\end{proofof}

\begin{proofof}{Theorem~\ref{th:FxMsol}}\label{pf:th:FxMsol}
    According to the results of Corollary~\ref{cor:bc}, we write $F(k)$ to mean $F(k-)=F(k+)$ and similarly for the derivative in $k$.
By defining 
    \begin{eqnarray}
        H_1 &=& (U^+_1 - U^-_1)^{-1} (e^{U^+_1 k} - e^{U^-_1 k}), \nonumber \\
        H_2 &=& M_0 (I - e^{U^-_1 k} + U^-_1 H_1), \nonumber \\
        H_3 &=& H_1 + \tD_1 H_2, \nonumber \\
        H_4 &=& - \lambda B_1 H_2, \nonumber 
    \end{eqnarray}
we can rewrite \eqref{eq:F1.exp} evaluated in $k$ as 
    \begin{eqnarray}
        F(k) 
        &=& F'(0) H_3 + \state{c-1} H_4 .
    \end{eqnarray}
By defining 
    \begin{eqnarray}
        H_5 &=& (U^+_1 - U^-_1)^{-1} (U^+_1 e^{U^+_1 k} - U^-_1 e^{U^-_1 k}), \nonumber \\
        H_6 &=& M_0 (U^-_1 e^{U^-_1 k} - U^-_1 H_5), \nonumber \\
        H_7 &=& H_5 - \tD_1 H_6, \nonumber \\
        H_8 &=& \lambda B_1 H_6, \nonumber 
    \end{eqnarray}
we can rewrite the derivative of Equation \eqref{eq:F1.exp} evaluated in $k$ as 
    \begin{eqnarray}\label{F'k-}
        F'(k) 
        &=& F'(0) H_7 + \state{c-1} H_8 . 
    \end{eqnarray}
By defining 
\begin{eqnarray}
    H_9 &=& (\tD_1-\tD_2) M_2 U^-_2 
    + \tD_1 (\tD_1-\tD_2) M_2, \nonumber
\end{eqnarray}
we can rewrite the derivative of Equation \eqref{eq:F2.exp} evaluated in $k$ as 
\begin{eqnarray}
    F'(k) 
    &=& F(k) U^-_2 - b_c \psic U^-_2 - \alpha_1 M_1 U^-_2  - \alpha_2 H_9 .
\end{eqnarray}
Then substituting the expression of $\alpha_2$ in \eqref{eq:alp2}, 
by employing also \eqref{eq:alp1} and \eqref{eq:alp0}, and defining
\begin{eqnarray}
    H_{10} &=&  U^-_2 - \lambda (I - B_2\tD_2^{-1}) H_9, \nonumber \\
    H_{11} &=&  M_1 U^-_2  + \tD_2^{-1} H_9, \nonumber \\
    H_{12} &=& H_{10} + \lambda (B_1\tD_1^{-1}\tD_2  - B_2) H_{11}, \nonumber \\
    H_{13} &=& \tD_2 H_{11}, \nonumber \\
    H_{14} &=& \lambda B_1  \tD_1^{-1} \tD_2 H_{11}, \nonumber 
\end{eqnarray}
we get 
\begin{eqnarray}\label{F'k+}
    F'(k) (I - H_9)
    &=& F(k) H_{12} - b_c \psic U^-_2 - F'(0) H_{13} + \state{c-1} H_{14} . 
\end{eqnarray}
By equating \eqref{F'k-} and \eqref{F'k+} and defining
\begin{eqnarray}
    H_{15} &=& U^-_2 (H_7 - H_7  H_9 - H_3 H_{12} + H_{13})^{-1},  \label{eq:H15} \\
    H_{16} &=& (H_{14} + H_4 H_{12} - H_8 + H_8 H_9)  (U^-_2)^{-1}  H_{15}, \label{eq:H16}
\end{eqnarray}
we get the first result in \eqref{eq:F'0}.
\begin{eqnarray}\label{eq:df0.2}
    F'(0) &=&  \state{c-1} H_{16} - b_c \psic   H_{15}
\end{eqnarray}
Taking the limit in \eqref{eq:F2.exp} and defining
\begin{eqnarray}
    H_{17} &=& \tD_2 M_1 - \lambda H_3 (B_1\tD_1^{-1}\tD_2  - B_2) M_1 , \nonumber \\
    H_{18} &=& \lambda  B_1 \tD_1^{-1} \tD_2 M_1 + \lambda H_4 (B_1\tD_1^{-1}\tD_2  - B_2) M_1 , \nonumber \\
    H_{19} &=& I -   H_{15} H_{17}, \label{eq:H19} \\
    H_{20} &=& H_{16} H_{17} - H_{18}, \label{eq:H20} 
\end{eqnarray}
we get the second result in \eqref{eq:Finfty}. 
\end{proofof}

\begin{proofof}{Theorem~\ref{th:bc.sol}}\label{pf:th:bc.sol}
    We would need the definition of the following rectangular matrices, for $0 \leq n \le c-1$:
    \begin{eqnarray*}
        \hat B_n(i,i) & = & (n-i+1)\mu_2, \;\;0 \le i \le n,\\
        \hat B_n(i,i-1) & = & i \mu_1, \;\;1 \le i \le n+1,\\
        \hat B_n(i,j) & = & 0 \;\;\;\mbox{for all other } (i,j),  \;\;0 \le j \le n .
    \end{eqnarray*}

    Using Equation \eqref{F.zero.cond} together with the balance equations \eqref{eq:piij0} and the normalization equation, we finally get 
\begin{eqnarray}
    \lambda \state{n} &=& \state{n+1} \hat B_0 ,\hspace{7.5em} n = 0, \label{eq.0}\\
    \state{n} (\lambda I + \mD_n ) &=& \state{n-1} \lambda \hat I + \state{n+1} \hat B_n ,\;\;\;  0 < n < c-1, \label{eq.n} \\
    \state{n} (\lambda I + \mD_n) &=& \state{n-1} \lambda  \hat I + F'(0) , \hspace{3.5em}  n = c-1, \label{eq.c-1}\\
    1 &=& F(\infty) \, \bone + \sum_{0 \le n \le c-1} \state{n} \, \bone ,\label{eq.norm}
\end{eqnarray}
where $F'(0)$ and $F(\infty)$  are given in \eqref{eq:F'0} and \eqref{eq:Finfty} respectively.
This system has $1 + (c + 1) c/2 $ equations and an equal number of unknowns.

Writing $\state{n}  = \state{n+1} \hat C_n$, we have, by \eqref{eq.0}, 
$ \hat C_0 =  \hat B_0 / \lambda $
and by \eqref{eq.n}, 
$ \hat C_n = \hat B_n (\lambda (I- \hat C_{n-1} \hat I) + \mD_n)^{-1}$, $0 < n < c-1$.
By \eqref{eq.c-1} we have
$\state{c-1}  = b_c \psic \hat C_{c-1}$
with 
$\hat C_{c-1}  = -   H_{15} (\lambda (I - \hat C_{c-2} \hat I) + \mD_{c-1}  - H_{16})^{-1}$ .

By defining 
\begin{equation}\label{eq:hat.H} 
    \hat H_{c-1}= \hat C_{c-1} , \;\;\; \hat H_{n} = \hat H_{n+1} \hat C_n , \;\;\;  0 < n < c-1,
\end{equation}
and using the normalization constraint \eqref{eq.norm} we get the result.
\end{proofof} 

\begin{proofof}{Lemma~\ref{lm:mean.computation}} 
\label{pf:lemma.meanVQT}
    By taking derivatives of Equations \eqref{eq:F1.exp} and \eqref{eq:F2.exp} we can compute the VQT density function as follows
    \begin{eqnarray}\label{eq:F1'.exp}
        F'(x) &=& (F'(0) + \alpha_0 M_0 U^-_1) (U^+_1 - U^-_1)^{-1} (U^+_1 e^{U^+_1 x} - U^-_1  e^{U^-_1 x}) \nonumber \\
        & & - \alpha_0 M_0 U^-_1 e^{U^-_1 x}, \;\;\; 0 \le x \le k , \\
        \label{eq:F2'.exp}
        F'(x) &=& (F(k) - b_c \psic - \alpha_1 M_1 - \alpha_2 (\tD_1-\tD_2) M_2) U^-_2 e^{U^-_2 (x-k)}  \nonumber \\
        & & - \alpha_2 \tD_1 \tQ_1(x-k) (\tD_1-\tD_2) M_2, \;\;  x \ge k .
    \end{eqnarray}
    We define the following matrix function
    \begin{equation}\label{eq:int.x.exp(Dx)}
        I(a,b;D) \defeq \int_a^b D x \; e^{D x} dx = (b e^{D b} - a e^{D a}) - D^{-1} (e^{D b}  - e^{D a}),
    \end{equation}
    that is 
    well defined on the set of non-singular matrices and that can be defined on the set of singular matrices by continuity. That is, if $\det(D)=0$, we set $I(a,b;D) = \lim_{t\to0} I(a,b;D+tI)$.
    
    Integrating the expressions \eqref{eq:F1'.exp} and \eqref{eq:F2'.exp} in their corresponding domains multiplied by $x$, we obtain the result in \eqref{eq:mean.computation} after summing up all components.
\end{proofof}

\bibliography{bib_waiting-dep-service}
\bibliographystyle{siam}
\end{document}

%% file: diagram.tex
\begin{tikzpicture}[
    scale=1.75,
    every label/.append style = {font = \scriptsize},
    every node/.append style = {font = \scriptsize},
    dot/.style = {minimum size = .1cm, inner sep = +0pt,
        shape = circle, draw = black, fill = black,
        label = {below: \({#1}\)}},
    larrow/.style = {->, shorten >= 2pt, shorten <= 2pt, transform canvas={yshift=-2pt}},
    rarrow/.style = {->, shorten >= 2pt, shorten <= 2pt, transform canvas={yshift=+2pt}},
    darrow/.style = {->, shorten >= 2pt, shorten <= 12pt},
    jump/.style = {draw=black, fill=white, inner sep=1pt,minimum size=5pt},
  ]
  
\ifx\c\undefined 
    \def\c{4}
\fi

\pgfmathparse{\c - 1}\let\n\pgfmathresult
\def\w{6}
\def\k{2}

\let\xlim\n
\foreach \x in {0,...,\xlim} {
  \pgfmathparse{\n - \x}\let\ylim\pgfmathresult
  \foreach \y in {0,...,\ylim} {
    \path (\x,\y) node[dot = {\x,\y}] (\x,\y) {};
  }
}

\pgfmathparse{\n - 1}\let\xlim\pgfmathresult
\foreach \x in {0,...,\xlim} {
  \pgfmathparse{\n - \x - 1}\let\ylim\pgfmathresult
  \foreach \y in {0,...,\ylim} {
    \path[larrow] (\x,\y) edge node[below=-2pt] {$\lambda$} ($ (\x,\y) + (1,0) $);
  }
}

\let\xlim\n
\foreach \x in {1,...,\xlim} {
  \pgfmathparse{\n - \x}\let\ylim\pgfmathresult
  \foreach \y in {0,...,\ylim} {
    \if\x1\def\s{$\mu_1$}\else\def\s{$\x\mu_1$}\fi
    \path[rarrow] (\x,\y) edge node[above=-3pt] {\s} ($ (\x,\y) - (1,0) $);
  }
}

\pgfmathparse{\n - 1}\let\xlim\pgfmathresult
\foreach \x in {0,...,\xlim} {
  \pgfmathparse{\n - \x}\let\ylim\pgfmathresult
  \foreach \y in {1,...,\ylim} {
    \if\y1\def\s{$\mu_2$}\else\def\s{$\y\mu_2$}\fi
    \path[darrow] (\x,\y) edge node[right=-3pt] {\s} ($ (\x,\y) - (0,1) $);
  }
}

\let\xlim\n
\foreach \x in {0,...,\xlim} {
  \pgfmathparse{\n - \x}\let\y\pgfmathresult
  \pgfmathparse{\w - 2}\let\dlim\pgfmathresult
  \draw[<-, shorten <= 2pt] (\x,\y) -- +(1,0);
  \foreach \d in {1,...,\dlim} {
    \draw[<-] ($(\x,\y)+(\d,0)$) -- +(1,0);
  }
  \draw[latex->] ($(\x,\y)+(\w,0)$) node [below right] {($w$,\x,\pgfmathprintnumber\y)} -- +(-1,0);
}

\draw[dashed] ($(0,\n) +(\k,0) +(-0.1,0.1)$) -- ($(\n,0) +(\k,0) +(0.1,-0.1)$) node [below right, font = \small] {$k$};

\draw[->] (0.2,-0.6) to[bend left] node[circle, draw=black, fill=white, inner sep=1pt, minimum size=5pt] {$i$} +(1,0) 
  node[right, font=\small] {$\lambda(i+1) \mu_1 \exp(-\Delta(i+1,c-1-i) x) \mathrm{d}x$ {\color{black} -- next exit of type $1$}};
\draw[->] (0.2,-0.9) to[bend left] node[draw=black, fill=white, inner sep=1pt, minimum size=5pt] {$i$} +(1,0) 
  node[right, font=\small] {$\lambda(c-1-i) \mu_2 \exp(-\Delta(i+1,c-1-i) x) \mathrm{d}x$ {\color{black} -- next exit of type $2$}};

\foreach \x in {0,...,\n} {
  \pgfmathparse{\n - \x}\let\y\pgfmathresult
  \pgfmathrnd\pgfmathroundto{\pgfmathresult}\let\s\pgfmathresult
  \pgfmathparse{greater(\s,0.5) * \k * 3/4 + greater(0.5,\s) * \w * 4 / 5}\let\dx\pgfmathresult
  \draw[->] (\x,\y) edge[bend left] node[circle,jump] {$\x$} +(\dx,0);
  \if\x\n\else
  \pgfmathparse{\dx + 1}\let\dx\pgfmathresult
  \draw[->] (\x,\y) edge[bend left] node[jump] {$\x$} +(\dx,-1);
  \fi
}

\foreach \x in {0,...,\n} {
  \pgfmathrnd\let\s\pgfmathresult
  \pgfmathparse{\x + \k * (1 + 6 * \s) / 8}\let\xs\pgfmathresult
  \pgfmathrnd\let\r\pgfmathresult
  \pgfmathparse{(\w - \k * (1 + 6 * \s) / 8) * (1 + 6 * \r) / 8}\let\dx\pgfmathresult
  \pgfmathparse{\n - \x}\let\y\pgfmathresult
  \draw[->] (\xs,\y) edge[bend left] node[circle,jump] {$\x$} +(\dx,0);
  \if\x\n\else
  \pgfmathparse{\dx + 1}\let\dx\pgfmathresult
  \draw[->] (\xs,\y) edge[bend left] node[jump] {$\x$} +(\dx,-1);
  \fi
}

\foreach \x in {0,...,\n} {
  \pgfmathrnd\let\s\pgfmathresult
  \pgfmathparse{\x + \k + (\w - \k) * (1 + 3 * \s) / 8}\let\xs\pgfmathresult
  \pgfmathrnd\let\r\pgfmathresult
  \pgfmathparse{(\w - \k - (\w - \k) * (1 + 3 * \s) / 8) * (3 + 4 * \r) / 8}\let\dx\pgfmathresult
  \pgfmathparse{\n - \x}\let\y\pgfmathresult
  \draw[->] (\xs,\y) edge[bend left] node[jump] {$\x$} +(\dx,0);
  \if\x0\else
  \pgfmathparse{\dx - 1}\let\dx\pgfmathresult
  \draw[->] (\xs,\y) edge[bend left] node[circle,jump] {$\x$} +(\dx,1);
  \fi
}

\end{tikzpicture}